\documentclass[10pt,twoside]{article}

\usepackage[paper=a4paper, left=3cm, right=3cm, top=2.5cm, bottom=3cm]{geometry}

\usepackage{amsmath,amssymb}

\usepackage{varioref}
\usepackage{hyperref}
\usepackage[noabbrev]{cleveref}

\usepackage{booktabs}

\usepackage{tabularx}
\newcolumntype{L}{>{\raggedright\arraybackslash}X}
\newcolumntype{R}{>{\raggedleft\arraybackslash}X}

\usepackage{etoolbox}\preto\subequations{\ifhmode\unskip\fi} 

\usepackage{tikz}
\usepackage{pgfplots}
\pgfplotsset{
/pgfplots/colormap={jet}{rgb255(0cm)=(0,0,128) rgb255(1cm)=(0,0,255) rgb255(3cm)=(0,255,255) rgb255(5cm)=(255,255,0) rgb255(7cm)=(255,0,0) rgb255(8cm)=(128,0,0)}}

\usepgfplotslibrary{external}
\tikzset{png export/.style={
   external/system call={
     lualatex \tikzexternalcheckshellescape
       -halt-on-error -interaction=batchmode -jobname "\image" "\texsource";
     convert -units pixelsperinch -density 600 "\image.pdf" "\image.png";
}}}
\tikzset{png export}
\usepgfplotslibrary{patchplots}
\usepackage{ifthen}
\newboolean{UseExternalPngs}
\setboolean{UseExternalPngs}{true}

\usepackage[ruled]{algorithm2e}
\usepackage{algorithmic}

\crefname{problem}{problem}{problems}
\crefname{algocf}{algorithm}{algorithms}

\DeclareMathOperator{\spn}{span} 

\definecolor{darkgreen}{rgb}{0.3,0.6,0.2}

\begin{document}

\title{POD model order reduction with space-adapted snapshots for incompressible flows}

\author{Carmen Gr\"a\ss{}le\thanks{Carmen Gr\"a\ss{}le, Department of Mathematics, Universit\"at Hamburg, Bundesstr. 55, 20146 Hamburg, Germany, carmen.graessle@uni-hamburg.de} \and Michael Hinze\thanks{Michael Hinze, Mathematisches Institut, Universit\"at Koblenz-Landau, Campus Koblenz, Universit\"atsstra\ss{}e 1, 56070 Koblenz, Germany, hinze@uni-koblenz.de} \and Jens Lang\thanks{Jens Lang, Department of Mathematics, Technische Universit\"at Darmstadt, Dolivostr. 15, 64293 Darmstadt, Germany, lang@mathematik.tu-darmstadt.de} \and Sebastian Ullmann\thanks{Graduate School of Computational Engineering, Technische Universit\"at Darmstadt, Dolivostr. 15, 64293 Darmstadt, Germany, ullmann@gsc.tu-darmstadt.de}}

\date{}

\newtheorem{problem}{Problem}
\newtheorem{remark}{Remark}

\maketitle

\begin{abstract}
We consider model order reduction based on proper orthogonal decomposition (POD) for unsteady incompressible Navier-Stokes problems, assuming that the snapshots are given by spatially adapted finite element solutions. We propose two approaches of deriving stable POD-Galerkin reduced-order models for this context. In the first approach, the pressure term and the continuity equation are eliminated by imposing a weak incompressibility constraint with respect to a pressure reference space. In the second approach, we derive an inf-sup stable velocity-pressure reduced-order model by enriching the velocity reduced space with supremizers computed on a velocity reference space. For problems with inhomogeneous Dirichlet conditions, we show how suitable lifting functions can be obtained from standard adaptive finite element computations. We provide a numerical comparison of the considered methods for a regularized lid-driven cavity problem.
\end{abstract}

\section{Introduction}

Many tasks in computational fluid dynamics involving incompressible and multi-phase flow are challenging since the underlying system of equations are expensive to solve. Two ways to decrease the associated simulation costs are spatially adaptive discretizations and model order reduction. Our approach to simulation-based model order reduction combines both strategies.\\
The novelty of this paper consists in the application of model order reduction based on proper orthogonal decomposition with space-adapted snapshots \cite{GH17,URL16} to the context of simulation of unsteady flow problems governed by the incompressible Navier-Stokes equations. As a result, the challenge arises of deriving a stable reduced-order model, since for space-adapted snapshots the weak divergence-free property only holds true in the respective adapted finite element space. For this reason, the contribution of this paper lies in proposing two approaches to formulating a stable reduced-order model. In the first approach, we use a projection of either the velocity snapshot data or the velocity POD basis onto a reference space. The projection is constructed in such a way that the resulting velocity POD modes are weakly divergence-free with respect to a pressure reference space. Consequently, the pressure term in the weak form of the Navier-Stokes system vanishes and the continuity equation is fulfilled by construction. This approach can be viewed as a generalization of the method of \cite{Sir87} to space-adapted snapshots. The second approach is a Galerkin projection of the primitive equations onto a POD space for the pressure field and an enriched POD space for the velocity field, in the spirit of \cite{BMQR15,RV07}. The enrichment functions are computed from the pressure POD in order to achieve inf-sup stability with respect to a pair of reference velocity and pressure spaces.

An efficient and sufficiently accurate reduction of the high-fidelity systems by POD reduced-order models based on space-adapted snapshots allows to use these models in a multi-query scenario like uncertainty quantification, where an ensemble of simulations is required to estimate statistical quantities, or optimal control, where a system of equations has to be solved repeatedly in order to find a minimum of a given cost functional. We intend to study this in future work.

Reduced-order modeling is applied to flow systems in the pioneering works \cite{BHL93,Pet89}. POD model order reduction for optimal control of fluids is studied in \cite{Rav00}, for example. An adaptive control strategy within optimal control of flows is given in \cite{AH01}. In order to adapt the reduced-order model for a flow control problem within the optimization, a trust-region POD framework is proposed in \cite{AFS00}. A theoretical investigation providing error estimates for POD approximations of a general equation in fluid dynamics is carried out in \cite{KV02}. Reduced basis methods using an offline/online procedure are applied to parametrized Navier-Stokes equations in \cite{QR07} addressing the pressure treatment and stability issues. POD-Galerkin reduced-order modeling for incompressible flows with stochastic Dirichlet boundary conditions is studied in \cite{Ull14}. We refer to \cite{LMQR14} for a general review on model order reduction for fluid dynamics in the case of static spatial discretizations. Recently, stabilization techniques for reduced-order methods for parametrized incompressible Navier-Stokes equations, where the full-order approximation is based on finite volume schemes, are investigated in \cite{SR18}.

A number of publications have considered model order reduction by projection onto a reduced space generated from space-adapted snapshots: Reduced basis methods with space-adapted snapshots are considered by \cite{ASU17,Ste2014} in the context of parametrized partial differential equations. The authors derive estimates of the error with respect to the infinite-dimensional truth solution. The key ingredient for these estimates is the use of a wavelet discretization scheme, which allows a numerical approximation of the dual norm of the infinite-dimensional residual. A different approach is taken by \cite{Yan16,Yan17,Yan18}, where bounds for the dual norm of the residual are provided for the case of minimum-residual mixed formulations of parametrized elliptic partial differential equations. An adaptive Galerkin finite element formulation is considered in \cite{URL16}, where computational issues of POD-Galerkin modeling in the presence of space-adapted snapshots are resolved by resorting to a common finite element mesh. Thus, an exact representation of the snapshots in the associated common finite element space is ensured. In the case of hierarchical, nested meshes, the construction of a common finite element mesh is given by an overlay of all adapted meshes and is cheap to construct. An a priori error analysis as well as an infinite-dimensional perspective in the context of evolution equations is provided by \cite{GH17}. This view also allows finite element discretizations in which the overlay of the adapted meshes leads to cut elements. For this case, a numerical implementation of the snapshot gramian is provided which, however, can be computationally demanding. In contrary to the just mentioned offline adaptive strategies, an online adaptive method is proposed in \cite{Car15} which provides a reduced-order analogon to $h$-refinement and is based on a splitting of the reduced basis vectors.

Our work is structured as follows: In \cref{sec:setting} we introduce the basic problem setting of an incompressible Navier-Stokes problem in strong and weak form. For ease of presentation, we first consider the setting with homogeneous Dirichlet boundary conditions. We provide an implicit Euler discretization in time and an adaptive Taylor-Hood finite element discretization in space. The adaptivity is achieved by combining residual-based error estimation, D\"orfler marking and newest vertex bisection. \Cref{sec:pod} introduces the fundamental concepts required to build a POD reduced basis from a set of functions. This section also defines an abstract reduced-order model, providing a framework for the following developments. \Cref{sec:velocityPOD} proposes a reduced-order model for the velocity field. It is based on a POD basis which is divergence-free in a weak sense with respect to a reference pressure space. A coupled velocity-pressure reduced-order model is introduced in \cref{sec:supremizers}. In order to ensure its stability, a set of supremizer functions is added to the velocity POD basis. A detailed presentation of the incorporation of inhomogeneous Dirichlet data is given in \cref{sec:inh_dirichlet}. Finally, the benchmark problem of a regularized lid-driven cavity flow serves as numerical test setting in \cref{sec:example}, in order to compare the methods regarding accuracy and computation time.

\section{Problem setting} \label{sec:setting}

\newcommand{\vel}[1][]{y^{#1}}  
\newcommand{\velx}{\vel_1}
\newcommand{\vely}{\vel_2}
\newcommand{\prs}[1][]{p^{#1}}  

\newcommand{\ti}[1][]{t^{#1}}
\newcommand{\velInit}{\vel_0}
\newcommand{\velDir}{y_\text{D}}

\newcommand{\hilbVel}{\mathcal V}
\newcommand{\hilbPrs}{\mathcal Q}
\newcommand{\divergence}{\nabla \cdot}

\newcommand{\inner}[4][]{#1(#2,#3#1)_{#4}}
\newcommand{\innerLtwo}[3][]{\inner[#1]{#2}{#3}{}}
\newcommand{\duality}[4][]{#1\langle#2,#3#1\rangle_{#4}}

\newcommand{\nSnapshots}{n}
\newcommand{\velh}[1][]{{\vel}_h^{#1}} 
\newcommand{\velhx}[1][]{{\vel}_{h,1}^{#1}}
\newcommand{\velhy}[1][]{{\vel}_{h,2}^{#1}}
\newcommand{\prsh}[1][]{{\prs}_h^{#1}} 
\newcommand{\snapVel}[1]{\vel[#1]}
\newcommand{\snapPrs}[1]{\prs[#1]}

\newcommand{\hilbSnapVel}[1][+]{V^{#1}} 
\newcommand{\hilbSnapPrs}[1][+]{Q^{#1}}

\newcommand{\tria}{\mathcal T}
\newcommand{\edge}{\mathcal E}
\newcommand{\triangulation}[2][h]{\mathcal T_{#1}^{#2}}
\newcommand{\initialtriangulation}{\mathcal T_h^\text{init}}

\newcommand{\normalder}{\nu_\Omega}

\newcommand{\testVel}{w}
\newcommand{\testPrs}{q}

\newcommand{\Hdiv}{H_\text{div}}

We consider an unsteady incompressible flow problem governed by the Navier-Stokes equations in a bounded domain $\Omega \subset \mathbb{R}^2$ with boundary $\partial\Omega$ over a time interval $[0,T]$ with $T>0$. The governing equations for the velocity field $\vel=(\velx,\vely)$ and pressure $\prs$ are
\begin{subequations}\label{NaSt} 
\begin{alignat}{4}
\vel_t + (\vel \cdot \nabla) \vel - Re^{-1} \Delta \vel + \nabla \prs &= f &\qquad& \text{in } (0,T) \times \Omega, \label{momentum}\\
\divergence \vel &= 0 &\qquad& \text{in } (0,T) \times \Omega, \label{continuity}\\
\vel &= 0      &\qquad& \text{in } (0,T)\times\partial\Omega,\\
\vel &= \velInit &\qquad& \text{in } \{0\}\times\Omega,
\end{alignat}
\end{subequations}
where $Re$ is the Reynolds number, $f$ denotes a given body force and $\velInit$ is an initial velocity field with $\divergence\velInit = 0$ in $\Omega$.

\subsection{Weak formulation}

The finite element and reduced-order models considered in this work are based on a weak form of the problem given by \eqref{NaSt}. We provide the necessary functional analytic framework by introducing the Hilbert spaces $\hilbVel=H_0^1(\Omega)$, $\hilbPrs = L_0^2(\Omega) = \{q \in L^2(\Omega): \int_\Omega q dx = 0\}$ and $W_\star^1(0,T;\hilbVel) = \{v \in L^2(0,T;\hilbVel): v_t \in L^1(0,T; \hilbVel')\}$. For clarity, we use the same notation for vector-valued functions, meaning that all components of a vector-valued function belong to the corresponding scalar function space. The same holds for vector- and scalar-valued operators. As short-hand notations we use $\innerLtwo{\cdot}{\cdot}:=\inner{\cdot}{\cdot}{L^2(\Omega)}$ and $\duality{\cdot}{\cdot}{} = \duality{\cdot}{\cdot}{\hilbVel',\hilbVel}$. We define $\inner{u}{v}{\hilbVel} = \innerLtwo{\nabla u}{\nabla v}$ and $\|v\|_\hilbVel^2=\inner{v}{v}{\hilbVel}$ for all $u,v\in\hilbVel$ and we set $\inner{p}{q}{\hilbPrs} = \innerLtwo{p}{q}$ and $\|q\|_\hilbPrs=\|q\|_{L^2(\Omega)}$ for all $p,q\in\hilbPrs$. We introduce the space $\Hdiv = \{\testVel\in L^2(\Omega)\colon \divergence \testVel = 0, (\testVel\cdot\normalder)|_{\partial\Omega} = 0\}$, where $\normalder$ denotes the outward unit normal vector. Finally, we introduce the notations $a(u,v) := Re^{-1}\inner{u}{v}{\hilbVel}$, $c(w,u,v) := \innerLtwo{(w \cdot \nabla) u}{v}$ and $b(v,q) := -\innerLtwo{q}{\divergence v}$.

The weak form of \eqref{NaSt} reads as follows: For given $f\in L^2(0,T;\hilbVel')$ and $\velInit\in\Hdiv$, find a velocity $\vel \in W_\star^1(0,T;\hilbVel)$ satisfying $\vel(0) = \velInit$ and a pressure $\prs \in L^2(0,T;\hilbPrs)$ such that
\begin{subequations}\label{NaSt_weak}
\begin{alignat}{4}
\frac{d}{dt} \innerLtwo{\vel (t)}{v} + c(\vel(t),\vel(t),v) + a(\vel(t),v) + b(v,\prs(t)) & = \duality{f(t)}{v}{} & \; \forall v\in \hilbVel, \label{momentum_weak}\\
b(\vel(t),q) & =  0 & \; \forall q \in \hilbPrs  \label{continuity_weak},
\end{alignat}
\end{subequations}
for almost all $t \in (0,T)$. For existence and uniqueness of a solution to \eqref{NaSt_weak} we refer to \cite[chapter 3, theorems 3.1 and 3.2]{Tem1979}.

\subsection{Discretization} \label{sec:timediscretization}

We first discretize in time and then discretize in space. This allows us to use a different finite element space at each time instance. We apply the implicit Euler scheme to discretize \eqref{NaSt_weak} in time. To this end, we introduce a time grid $0=\ti[0] < \dots < \ti[\nSnapshots] = T$ with $\nSnapshots \in \mathbb{N}$. For simplicity, we assume an equidistant spacing with a time step size $\Delta t = T/\nSnapshots$. The time-discrete system consists in finding sequences $\vel[1],\dots,\vel[\nSnapshots]\in\hilbVel$ and $\prs[1],\dots,\prs[\nSnapshots]\in\hilbPrs$, for given $\vel[0] = \velInit \in \Hdiv$, satisfying the system
\begin{subequations}\label{NaSt_weakTimeDiscrete}
\begin{alignat}{4}
 \innerLtwo[\Big]{\frac{\vel[j]-\vel[j-1]}{\Delta t}}{v} + c(\vel[j],\vel[j],v) + a(\vel[j],v) + b(v,\prs[j]) & =  \duality{f(t^j)}{v}{} &\quad& \forall v\in \hilbVel,  \label{momentum_weakTimeDiscrete}\\
b(\vel[j],q) & =  0&\quad& \forall q \in \hilbPrs \label{continuity_weakTimeDiscrete}
\end{alignat}
\end{subequations}
for $j=1,\dots,n$. Note that we have applied the box rule in order to approximate the right-hand side time integral. An initial pressure field can be obtained from an additional pressure Poisson equation, if required, see e.g. \cite{HeywoodRannacher1982}.

For the spatial discretization, we utilize adaptive finite elements based on LBB stable $\mathbb{P}_2-\mathbb{P}_1$ Taylor-Hood elements. For each time instance, we use spatially adapted finite element spaces $\{\hilbSnapVel[1], \dots, \hilbSnapVel[\nSnapshots]\} \subset \hilbVel$ and $\{\hilbSnapPrs[1], \dots, \hilbSnapPrs[\nSnapshots]\} \subset \hilbPrs$, which we specify in \cref{sec:adaptive}. The fully discrete Navier-Stokes problems read as follows: Given $\velh[0] = \velInit\in\Hdiv$, find $\velh[1]\in\hilbSnapVel[1],\dots,\velh[\nSnapshots]\in\hilbSnapVel[\nSnapshots]$ and $\prsh[1]\in \hilbSnapPrs[1],\dots,\prsh[\nSnapshots]\in \hilbSnapPrs[\nSnapshots]$ such that
\begin{subequations}\label{NaSt_fullydisc} 
\begin{alignat}{4}
  \innerLtwo[\Big]{\frac{\velh[j]-\velh[j-1]}{\Delta t}}{v} + c(\velh[j],\velh[j],v) + a(\velh[j],v) + b(v,\prsh[j]) & =   \duality{f(t^j)}{v}{}  \label{momentum_fullydisc} &\quad& \forall v \in \hilbSnapVel[j],\\
b(\velh[j],q) & =  0 \label{continuity_fullydisc} &\quad& \forall q \in \hilbSnapPrs[j]
\end{alignat}
\end{subequations}
for $j=1,\dots,\nSnapshots$. In each step of the implicit Euler method we compute an inner product of the velocity $\velh[j-1]$ at the previous time level with test functions $v \in \hilbSnapVel[j]$ of the current time level. Alternatively, it is possible to interpret the inner product as an $L^2(\Omega)$-projection of $\velh[j-1]$ onto $\hilbSnapVel[j]$ under a weak divergence-free constraint with respect to $\hilbSnapPrs[j]$, see \cite[Lemma 4.1]{BesierWollner2012}. For existence of a unique solution to \eqref{NaSt_fullydisc}, we refer to \cite[Chapter 3, \textsection 5 Scheme 5.1]{Tem1979}.

\subsection{Adaptive finite element method} \label{sec:adaptive}

In the following, we describe the choice of the mixed finite element pairs $(\hilbSnapVel[1],\hilbSnapPrs[1]),$ $\dots,(\hilbSnapVel[\nSnapshots],\hilbSnapPrs[\nSnapshots])$. As a starting point we define an initial finite element grid $\initialtriangulation$. We obtain adapted grids $\triangulation{j}$ by refining this initial grid. For each adapted grid, we can define a corresponding Taylor-Hood finite element pair $(\hilbSnapVel[j],\hilbSnapPrs[j])$. The procedure that leads to the individual finite element pairs for a given initial grid can be described by the standard \emph{solve-estimate-mark-refine} cycle. The details for each of these steps are provided in \cref{Alg:AdFE}, an explanation is given below.

\begin{algorithm}[htbp]
\caption{Adaptive finite element algorithm.}
\label{Alg:AdFE}
\begin{algorithmic}[1]
\renewcommand{\algorithmicrequire}{\textbf{Input:}}
\renewcommand{\algorithmicensure}{\textbf{Output:}}
\REQUIRE Initial mesh $\initialtriangulation$, tolerance $\varepsilon > 0$, number of time instances $n$, refinement parameter $\theta\in(0,1)$. 
\ENSURE $\velh[1]\in\hilbSnapVel[1],\dots,\velh[\nSnapshots]\in\hilbSnapVel[\nSnapshots]$ and $\prsh[1]\in\hilbSnapPrs[1],\dots,\prsh[\nSnapshots]\in\hilbSnapPrs[\nSnapshots]$
\STATE Set $\triangulation{1} := \initialtriangulation$.
\FOR{$j=1, \dots, \nSnapshots$} 
  \LOOP
    \STATE Define $(\hilbSnapVel[j],\hilbSnapPrs[j])$ as the Taylor-Hood finite element pair corresponding to $\triangulation{j}$.
    \STATE \emph{Solve} \eqref{NaSt_fullydisc} for given $j$.
    \STATE \emph{Estimate} the error contributions $\eta_\tria^j$ from \eqref{error_contributions} for given $j$ for all $\tria \in \triangulation{j}$.
    \IF{$\sum_{\tria \in \triangulation{j}} \eta_{\tria}^j < \varepsilon$} 
    \STATE Construct $\triangulation{j+1}$ by coarsening $\triangulation{j}$ once and uniting with $\initialtriangulation$.
      \STATE \textbf{break} 
    \ENDIF
    \STATE \emph{Mark} the smallest set $ \triangulation[hD]{j} \subset \triangulation{j}$ which fulfills the D\"orfler criterion \eqref{eq:doerfler}.
    \STATE \emph{Refine} $\triangulation{j}$ using newest vertex bisection for given $\triangulation[hD]{j}$.
  \ENDLOOP
\ENDFOR
\end{algorithmic}
\end{algorithm}

For each $j$, the first part of the adaptive procedure is the \emph{solution} of the system of equations \eqref{NaSt_fullydisc} given a mixed finite element pair $(\hilbSnapVel[j],\hilbSnapPrs[j])$. The solution of \eqref{NaSt_fullydisc} leads to a non-linear algebraic saddle point problem, which must be solved for the velocity and pressure at the new time instance, given the velocity at the old time instance. We solve the non-linear system with Newton's method, using a standard sparse direct solver for the solution of the linear systems in each Newton iteration.

The error \emph{estimation} relies on a residual based a posteriori error estimator in the spirit of \cite{AO2000}. In particular, we obtain error indicators for the spatial error at each time step by adding the discrete time derivative to an error estimator for the stationary Navier-Stokes problem in \cite[section 4.4]{John2016}, or, equivalently, adding a convection term to an error estimator for the unsteady Stokes problem in \cite[section 5.4]{V2013}. The resulting estimator can also be found in \cite[section IV.2.2.]{V18}. It is given by
\begin{align}\label{error_contributions}
  \eta_\tria^j &=
  \Bigg(
      h_\tria^2\Bigg\|\displaystyle\frac{\velh[j]-\velh[j-1]}{\Delta t}
    + \velh[j] \cdot \nabla \velh[j]
    - Re^{-1} \Delta \velh[j]
    + \nabla \prsh[j]
    - f(t^j) \Bigg\|_{L^2(\tria)}^2\notag\\
    &  +\Big\|\nabla\cdot\velh[j]\Big \|_{L^2(\tria)}^2 
    +\frac{1}{2}\sum_{\edge\in\partial\tria\setminus\partial\Omega}\!\!\!h_{\edge}\Big\|\big[
    - Re^{-1}\nabla \velh[j] \cdot \normalder
    + \prsh[j] \normalder
  \big]_\edge\Big\|_{L^2(\edge)}^2
  \Bigg)^\frac{1}{2}
\end{align}
for all $\tria\in\triangulation{j}$ and for $j=1,\dots,\nSnapshots$, assuming $f(t^j) \in L^2(\Omega)$. Here, $h_\tria^2$ is the triangle area, $h_{\edge}$ is the edge length and $[\cdot]_\edge$ denotes a jump over the edge $\edge$.
 
We use the D\"orfler criterion \cite{D1996} as a \emph{marking} strategy. This means, for refinement we select the smallest subset $\triangulation[hD]{j}$ of $\triangulation{j}$ fulfilling the requirement
\begin{equation}\label{eq:doerfler}
    \sum_{\tria\in\triangulation[hD]{j}}\eta_{\tria}^j\geq(1-\theta)\sum_{\tria\in \triangulation{j}} \eta_{\tria}^j,
\end{equation}
where $\theta \in (0,1)$ is the refinement parameter.

As a \emph{refinement} procedure, we use the newest vertex bisection method \cite{Mit1991}, which has the advantage that the resulting meshes are nested. The smallest common mesh of two adapted meshes is their overlay \cite{CKNS2008,Ste2007}.

In order to construct the starting mesh for the next time step, we mark every triangle of the current triangulation $\triangulation{j}$ and perform one coarsening step. We unite the result with $\initialtriangulation$ in order to guarantee that the mesh does not become coarser than the initial mesh. This strategy ensures that it is possible to reach the initial mesh $\initialtriangulation$ after a finite number of time steps. Although this choice might lead to a finer triangulation compared to starting from the initial triangulation $\initialtriangulation$ in each time step, we expect the advantage that only a small number of refinement steps are needed when proceeding from one time step to the other.

\section{POD-Galerkin modeling} \label{sec:pod}

\newcommand{\velPOD}[1][]{{\vel}_\nPod^{#1}}
\newcommand{\prsPOD}[1][]{{\prs}_\nPod^{#1}}
\newcommand{\PODcoefvel}{\eta}
 
\newcommand{\hilbSnap}[1][]{X_{#1}}

\newcommand{\snap}[1][]{u^{#1}} 
\renewcommand{\pod}[1][]{\phi_{#1}}

\newcommand{\nPod}{R}
\newcommand{\iPod}{r}

\newcommand{\nPodVel}{\nPod_{\vel}}
\newcommand{\nPodPrs}{\nPod_{\prs}}

\newcommand{\hilbPodVel}[1][\nPod]{V_{#1}}
\newcommand{\hilbPodPrs}[1][\nPod]{Q_{#1}}

\newcommand{\podVel}[1]{\phi^{#1}}
\newcommand{\podPrs}[1]{\psi^{#1}}

\newcommand{\weight}[2][]{\alpha_{#1}^{#2}}

\newcommand{\dirVelh}[1][]{\hat\vel_h^{#1}}

\newcommand{\stokesVel}[1]{\hat g^{#1}}
\newcommand{\homVel}[1]{\hat y^{#1}}

\newcommand{\hilbRefVel}{\tilde V}
\newcommand{\hilbRefPrs}{\tilde Q}

In practice, solving \eqref{NaSt_fullydisc} can easily lead to large non-linear algebraic systems of equations, which are computationally expensive to solve. For this reason, we apply model order reduction in order to replace the high-dimensional systems of equations by a low-dimensional approximation, which represents the original problem reasonably well. We use proper orthogonal decomposition (POD) in order to provide low-dimensional approximation spaces and use them in a Galerkin framework to derive the reduced-order models. In the following, we first provide a general description of POD. Then, we introduce an abstract Galerkin model, which provides a common foundation for the concrete models described in \cref{sec:velocityPOD,sec:supremizers}.

In order to formulate the POD, assume a set of functions $\snap[1] , \dots, $ $\snap[\nSnapshots] \in\hilbSnap$ is given, where $\hilbSnap$ is a Hilbert space. In the context of model order reduction, these functions are usually called snapshots. They could, for instance, be infinite-dimensional velocity or pressure fields of the time-discrete problem \eqref{NaSt_weakTimeDiscrete}, or corresponding finite-dimensional approximations. In principle, the number of time instances and the number of snapshots could be chosen differently. However, for the sake of simplicity we take the same number $n$ of snapshots as the number of time instances in the scope of this work.

The POD method consists of finding functions $\pod[1],\dots,\pod[\nPod]\in\hilbSnap$ with $\nPod\leq\nSnapshots$, which solve the equality constrained minimization problem
\begin{equation}\label{minPOD}
\underset{\pod[1] , ... , \pod[\nPod]}{\text{min }}
 \displaystyle\sum_{j=1}^\nSnapshots \weight{j} \left\| \snap[j] - \displaystyle\sum_{i=1}^\nPod 
 \inner{\snap[j]}{\pod[i]}{\hilbSnap} \; \pod[i] 
 \right\|_{\hilbSnap}^2  \text{ s.t. }  \inner{\pod[i]}{\pod[j]}{\hilbSnap} = 
 \delta_{ij} \text{ for } 1 \leq i, j \leq \nPod,
 \end{equation}
with $\{\weight{j}\}_{j=1}^n$ denoting non-negative weights and $\delta_{ij}$ the Kronecker symbol. This minimization problem can be solved using a generalized eigenvalue decomposition of a snapshot Gramian, see \cite{Sir87} for instance. The functions $\pod[1],\dots,\pod[\nPod]$ are called POD basis functions. If $\nPod\leq\dim(\spn(\snap[1] , \dots, \snap[\nSnapshots]))$, then $\pod[1],\dots,\pod[\nPod]\in \spn(\snap[1] , \dots, \snap[\nSnapshots])$. An $\nPod$-dimensional POD space is given by $\spn(\pod[1],\dots,\pod[\nPod])$.

In order to derive a reduced-order model (ROM) of the time-discrete weak form \eqref{NaSt_weakTimeDiscrete}, we introduce abstract reduced spaces $\hilbPodVel\subset\hilbVel$ and $\hilbPodPrs\subset\hilbPrs$ for the velocity and pressure, respectively. The index $\nPod$ is related, but not necessarily equal, to the dimensions of $\hilbPodVel$ or $\hilbPodPrs$. Concrete choices of $\hilbPodVel$ and $\hilbPodPrs$ are provided in \cref{sec:velocityPOD,sec:supremizers}.

Replacing the original spaces $\hilbVel$ and $\hilbPrs$ in \eqref{NaSt_weakTimeDiscrete} with the respective reduced spaces $\hilbPodVel$ and $\hilbPodPrs$ leads to the following  abstract reduced-order problem: For given $\velPOD[0]=\velInit\in\Hdiv$, find $\velPOD[1],\dots,\velPOD[\nSnapshots]\in\hilbPodVel$  and $\prsPOD[1],\dots,\prsPOD[\nSnapshots]\in\hilbPodPrs$ such that
\begin{subequations}\label{tmp}
\begin{alignat}{4}
  \innerLtwo[\Big]{\frac{\velPOD[j]-\velPOD[j-1]}{\Delta t}}{v} + c(\velPOD[j],\velPOD[j],v) + a(\velPOD[j],v) + b(v,\prsPOD[j])& =   \duality{f(t^j)}{v}{}  &\quad& \forall v \in \hilbPodVel,\label{tmp_a}\\
b(\velPOD[j],q) & =  0&\quad& \forall q \in \hilbPodPrs\label{tmp_b}
\end{alignat}
\end{subequations}
for $j=1,\dots,\nSnapshots$. We note that \eqref{tmp} constitutes a system of algebraic equations for the expansion coefficients of the reduced solutions. In this view, the system does not depend on the full spatial dimension of \eqref{NaSt_fullydisc}, see \cite[part III, section 1]{Sir87}.

The stability of \cref{tmp} is not guaranteed for all pairs of $\hilbPodVel$ and $\hilbPodPrs$. In the following, we provide two choices for $(\hilbPodVel,\hilbPodPrs)$ which result in stable reduced-order models.

The first approach is a velocity reduced-order model, presented in \cref{sec:velocityPOD}. It relies on a reference pressure finite element space paired with a velocity POD space, where the POD basis functions are weakly divergence-free with respect to the reference pressure space. This enables a cancellation of the pressure term in \eqref{tmp_a} and the continuity equation \eqref{tmp_b} is fulfilled by construction. The stability of the resulting system is then given by \cite{CF88,Tem1979}.

The second approach is a velocity-pressure reduced-order model, investigated in \cref{sec:supremizers}. It combines a pressure POD space with a velocity POD space which is  augmented by supremizer functions in order to achieve stability, see \cite{BMQR15,RV07}.

In the following, we make use of a reference velocity space $\hilbRefVel \subset \hilbVel$ and an associated reference pressure space $\hilbRefPrs \subset \hilbPrs$, such that the pair $(\hilbRefVel,\hilbRefPrs)$ is inf-sup stable. We do not impose further assumptions on these reference spaces and they are therefore kept general. One choice could be to take the common finite element spaces $\hilbRefVel$ and $\hilbRefPrs$ which contain all finite element spaces, i.e. $\hilbSnapVel[1], \dots, \hilbSnapVel[\nSnapshots]\subset\hilbRefVel$ and $\hilbSnapPrs[1],\dots,\hilbSnapPrs[\nSnapshots]\subset\hilbRefPrs$. For example, if the meshes $\mathcal{T}_h^1, \dots, \mathcal{T}_h^\nSnapshots$ are obtained by successive newest vertex bisections applied to a common coarse grid, then the spaces $\hilbRefVel$ and $\hilbRefPrs$ can be constructed using the overlay of all meshes. However, it is also possible to choose $\hilbRefVel$ and $\hilbRefPrs$ which are independent of the snapshot spaces.

\section{Velocity reduced-order model} \label{sec:velocityPOD}

\newcommand{\snapVelMod}[1]{\hat y^{#1}}
\newcommand{\snapVelProj}[1]{\tilde y^{#1}}
\newcommand{\snapPrsMod}[1]{\hat p^{#1}}

\newcommand{\podVelMod}[1]{\hat \phi^{#1}}

\newcommand{\hilbRefVelDiv}{\hilbRefVel_\text{div}}

\newcommand{\dir}[1][]{g^{#1}}
\newcommand{\dirh}[1][]{g^{#1}_h}
\newcommand{\dirDiv}[1][]{\hat g^{#1}}

\newcommand{\dirSnap}{\hat u}
\newcommand{\lerSnap}{\tilde u}
\newcommand{\lerayInh}[1][]{P_{#1}}
\newcommand{\lagrange}[1]{\mathbb I_{#1}}

Our goal is to derive a reduced-order model which only contains the velocity as an unknown. This can be achieved with a weakly divergence-free POD basis. If the given snapshots are weakly divergence-free with respect to one and the same test space, then this property carries over to the POD basis, see \cite{IR98,Sir87,VP05}. An analysis of a general POD-Galerkin flow model can be found in \cite{KV02}, for example.

The challenge in the context of adaptive spatial discretizations is that ``weakly divergence-free'' refers to the test space, which can be a different space at each time instance, in general. In particular, the solutions $\velh[j]\in\hilbSnapVel[j], j=1, \dots, \nSnapshots$ of \eqref{NaSt_fullydisc} fulfill a weak divergence-free property with respect to the corresponding pressure spaces $\hilbSnapPrs[j]$:
\[
  b(\velh[j],\testPrs) = 0\qquad\forall \testPrs\in\hilbSnapPrs[j],\;j=1,\dots,\nSnapshots.
\]
However, $\velh[j]$ is not necessarily weakly divergence-free with respect to the other spaces $\hilbSnapPrs[i]$ for $i\neq j$. As a consequence, no weak divergence-free property can be guaranteed for arbitrary linear combinations of snapshots. This means, if we compute a POD of $\velh[1],\dots,\velh[\nSnapshots]$, then the resulting POD basis functions are not necessarily weakly divergence-free. Therefore, we use $\velh[1],\dots,\velh[\nSnapshots]$ to construct a \emph{modified} velocity POD basis which is weakly divergence-free with respect to a reference pressure space. This allows an elimination of the pressure term in \eqref{tmp_a}, and the continuity equation \eqref{tmp_b} is fulfilled by construction. As a result, we directly obtain a velocity reduced-order model. 

We introduce two approaches to constructing a suitable modified velocity POD basis. The first approach is based on projected snapshots, meaning that we first project the snapshots such that they are weakly divergence-free and then compute a POD basis from the projected snapshots (\textit{first-project-then-reduce}, \cref{sec:approach_i}). The second approach is based on projected POD basis functions, implying that we first compute a POD basis from the original snapshots and then project the POD basis functions such that they fulfill a weak divergence-free property (\textit{first-reduce-then-project}, \cref{sec:approach_ii}).

\subsection{Optimal projection onto a weakly divergence-free space} \label{sec:leray}

Our aim is to define weakly divergence-free approximations in a more general sense. We provide a procedure that can be applied to problems with inhomogeneous Dirichlet data as well. To this end, we introduce a Dirichlet lifting function $\dir$, which will be specified in \cref{sec:inh_dirichlet}. In the context of homogeneous Dirichlet data, we set $\dir=0$. We would like to approximate a function $\snap\in\hilbSnap$ with a function $\lerSnap \in \hilbRefVel\subset \hilbSnap$ such that $\lerSnap + \dir$ is weakly divergence-free with respect to a space $\hilbRefPrs$. More precisely, we want to solve the following equality constrained minimization problem:
\begin{problem} \label{problem:lerayMinimization}
For given $\snap \in \hilbSnap$ and sufficiently smooth $\dir$, find $\lerSnap\in\hilbRefVel\subset \hilbSnap$ which solves
\[
  \min_{v\in\hilbRefVel}\frac12\|v-\snap\|_{\hilbSnap}^2\quad\text{subject to}\quad b(v+\dir,\testPrs)=0\quad\;\forall \testPrs\in\hilbRefPrs.
\]
\end{problem}
Note that \cref{problem:lerayMinimization} has a unique solution $\lerSnap\in\hilbRefVel$. We define the projection $\snap \mapsto \lerayInh[\dir](\snap) = \lerSnap$. In order to compute the solution to \cref{problem:lerayMinimization}, the usual Lagrange approach can be followed, see e.g. \cite{hpu}. The resulting system is the following saddle point problem:
\begin{problem} \label{problem:leraySaddlePoint}
For given $\snap \in \hilbSnap$ and sufficiently smooth $\dir$, find $\lerSnap \in \hilbRefVel\subset \hilbSnap$ and $\lambda\in\hilbRefPrs$ such that
\begin{alignat*}{2}
\inner{\lerSnap}{\testVel}{\hilbSnap} + b(\testVel,\lambda) &= \inner{\snap}{\testVel}{\hilbSnap} &\qquad& \forall\testVel\in\hilbRefVel,\\
b(\lerSnap,\testPrs) &= -b(\dir,\testPrs) && \forall\testPrs\in \hilbRefPrs.
\end{alignat*}
\end{problem}
Note that the solution $\lerSnap \in \hilbRefVel$ to \cref{problem:leraySaddlePoint} is  unique. If $\hilbRefVel$ and $\hilbRefPrs$ form an inf-sup stable pair of spaces, then the uniqueness of $\lambda\in\hilbRefPrs$ in \cref{problem:leraySaddlePoint} is given. If we choose $\hilbSnap=\hilbVel$, then \cref{problem:leraySaddlePoint} is a Stokes problem. For the choice $\hilbSnap=L^2(\Omega)$, \cref{problem:leraySaddlePoint} is a weak formulation of the Leray projection, see e.g.\ \cite{CF88,Tem1979}.

\subsection{Reduced-order modeling based on projected snapshots}\label{sec:approach_i}

The basic idea is to project the original velocity solutions of \eqref{NaSt_fullydisc} in order to obtain functions which are weakly divergence-free with respect to the reference pressure space $\hilbRefPrs$. Consequently, a resulting POD basis inherits this property by construction (\textit{first-project-then-reduce}). For given snapshots $\velh[1] \in \hilbSnapVel[1], \dots, \velh[\nSnapshots] \in \hilbSnapVel[\nSnapshots]$, we solve \cref{problem:leraySaddlePoint} with $\dir=0$. Then, the projected snapshots $\snapVelProj{1} = \lerayInh[0](\velh[1]), \dots, \snapVelProj{\nSnapshots} = \lerayInh[0](\velh[\nSnapshots])$ live in
\begin{equation}
  \hilbRefVelDiv:=\{v \in \hilbRefVel: b(v,\testPrs)=0 \quad \forall \testPrs \in \hilbRefPrs\}.\label{eq:hilbRefVelDiv}
\end{equation}
From these projected snapshots $\snapVelProj{1}, \dots, \snapVelProj{\nSnapshots}$, we compute a POD basis according to \eqref{minPOD} with Hilbert space $\hilbSnap=\hilbVel$, snapshot weights $\weight{1}=\dots=\weight{\nSnapshots}=\Delta t$ and POD dimension $\nPodVel  \leq \dim (\spn (\snapVelProj{1}, \dots, \snapVelProj{\nSnapshots}))$. The resulting POD space $\hilbPodVel := \spn(\podVel{1}, \dots, \podVel{\nPodVel})$ fulfills the property $\hilbPodVel \subset \hilbRefVelDiv \subset \hilbRefVel$. Thus, $b(\testVel,\testPrs) = 0$ holds true for all $\testVel \in \hilbPodVel$ and all $\testPrs \in \hilbRefPrs$. Consequently, for this choice of $\hilbPodVel$ in \eqref{tmp}, together with $\hilbPodPrs=\hilbRefPrs$, the pressure term vanishes and the continuity equation is fulfilled by construction. The resulting velocity ROM reads as follows: For given $\velPOD[0]=\velInit\in\Hdiv$, find $\velPOD[1], \dots, \velPOD[\nSnapshots] \in \hilbPodVel$ such that
\begin{equation}\label{NaSt_velPod_fullyreduced} 
  \innerLtwo[\Big]{\frac{\velPOD[j]-\velPOD[j-1]}{\Delta t}}{\testVel} + c(\velPOD[j],\velPOD[j],\testVel) + a(\velPOD[j],\testVel)  =  \duality{f(t^j)}{v}{}  \quad \forall \testVel \in \hilbPodVel,
\end{equation}
 for $j=1,\dots,\nSnapshots$. 

Concerning the computational complexity, we note that \cref{problem:lerayMinimization} has to be solved for each snapshot, which means the solution of a saddle point problem with reference spaces $\hilbRefVel$ and $\hilbRefPrs$, followed by the computation of a POD basis for the velocity field. Concerning the online computational costs, in each time step we solve a non-linear algebraic system of equations with Newton's method. In each Newton step, we need to build a Jacobian matrix and a right-hand side and, subsequently, solve a dense linear system using a direct method. Therefore, we find that solving the velocity ROM \eqref{NaSt_velPod_fullyreduced} is of order $\mathcal{O}(\nPodVel^3)$.

\subsection{Reduced-order modeling based on projected POD basis functions}\label{sec:approach_ii}

The basic idea of this approach is to compute a POD basis from the original velocity solutions of \eqref{NaSt_fullydisc} and project the resulting POD basis functions onto a weakly divergence-free space (\textit{first-reduce-then-project}). For given snapshots $\velh[1] \in \hilbSnapVel[1], \dots, \velh[\nSnapshots] \in \hilbSnapVel[\nSnapshots]$, we define interpolated snapshots $\snapVelMod{j} = \lagrange{\hilbRefVel}\velh[j]\in\hilbRefVel$ for $j=1, \dots, \nSnapshots$, where $\lagrange{\hilbRefVel}$ denotes the Lagrange interpolation operator onto the velocity reference space $\hilbRefVel$. From these interpolated snapshots, we compute POD basis functions $\podVelMod{1}, \dots, \podVelMod{\nPodVel} \in \hilbRefVel$ according to \eqref{minPOD} with Hilbert space $\hilbSnap=\hilbVel$, snapshot weights $\weight{1}=\dots=\weight{\nSnapshots}=\Delta t$ and POD dimension $\nPodVel  \leq \dim (\spn (\snapVelMod{1}, \dots, \snapVelMod{\nSnapshots}))$. These POD basis functions in general do not live in $\hilbRefVelDiv$. Thus, they are projected onto the space $\hilbRefVelDiv$ by solving \cref{problem:leraySaddlePoint} with $\dir=0$. Then, the projected POD basis functions $\podVel{1} = \lerayInh[0](\podVelMod{1}), \dots, \podVel{\nPodVel} = \lerayInh[0](\podVelMod{\nPodVel})$ live in $\hilbRefVelDiv$. Choosing $\hilbPodVel:=\spn(\podVel{1}, \dots, \podVel{\nPodVel})$ in \eqref{tmp} together with $\hilbPodPrs=\hilbRefPrs$ leads to a velocity ROM of the form \eqref{NaSt_velPod_fullyreduced}. Note that, in general, the reduced space constructed in this approach does not coincide with the reduced space constructed according to \cref{sec:approach_i}.

The computational complexity of the approach described in this subsection comprises the computation of a POD basis and, afterwards, the solution of \cref{problem:lerayMinimization} for each POD basis function, i.e. $\nPodVel$ times. This makes the current approach cheaper than the approach of \cref{sec:approach_i}, which required $\nSnapshots$ solutions of \cref{problem:lerayMinimization}, and $\nPodVel\leq\nSnapshots$. Otherwise the costs of setting up and solving the reduced-order model are equivalent.

\begin{remark}
Obviously, the velocity ROM \eqref{NaSt_velPod_fullyreduced} only depends on the velocity variable and the pressure is eliminated. However, many applications require an approximate pressure field. For example, the pressure is needed for the computation of lift and drag coefficients of an airfoil or for low-order modeling of shear flows, see e.g.\ \cite{NPM05}. In the case of snapshot generation on static spatial meshes, it is possible to reconstruct a reduced pressure afterwards by solving a discrete reduced pressure Poisson equation, see \cite{Ull14}, for example. A transfer of this concept to the case of space-adapted snapshots is not carried out within the scope of this work. In the following \cref{sec:supremizers}, we introduce a reduced-order model which depends on both velocity and pressure. Thus, it delivers directly a reduced pressure approximation without any post-processing recovery.
\end{remark}

\section{Velocity-pressure reduced-order model} \label{sec:supremizers}

\newcommand{\hilbPodSupVel}[1][\nPod]{V_{#1}}
\newcommand{\hilbPodSupPrs}[1][\nPod]{Q_{#1}}

\newcommand{\podSup}[1]{\bar\phi^{#1}}

\newcommand{\hilbRefVelDual}{\hilbRefVel'}

\newcommand{\infsupCom}{\tilde \beta_h}
\newcommand{\infsupRom}{\beta_R}
 
\newcommand{\bilDivergence}[2]{b(#1,#2)}

\newcommand{\supMap}{\mathbb T}
\newcommand{\infstack}[2]{\underset{\substack{\vphantom{\hilbPodPrs}{#1}\\\mathstrut{#2}}}{\mathstrut\inf}}
\newcommand{\supstack}[2]{\underset{\substack{\vphantom{\hilbPodPrs}{#1}\\\mathstrut{#2}}}{\mathstrut\sup}}

We want to derive a POD-Galerkin model which can be solved for reduced-order representations of the velocity and pressure fields. To this end, we need to choose reduced spaces such that the reduced-order model resulting from \eqref{tmp} becomes inf-sup stable. We obtain a suitable pressure reduced space by a truncated POD of a set of pressure snapshots. For the velocity reduced space, we take a truncated POD basis of a set of velocity snapshots and add suitable functions to guarantee the fulfillment of the inf-sup stability criterion, compare \cite{BMQR15,RV07}.

\subsection{POD spaces} \label{sec:supremizersPod}

Let $\velh[1]\in\hilbSnapVel[1],\dots,\velh[\nSnapshots]\in\hilbSnapVel[\nSnapshots]$ and $\prsh[1]\in \hilbSnapPrs[1],\dots,\prsh[\nSnapshots]\in \hilbSnapPrs[\nSnapshots]$ be the solutions of the fully discrete problem \eqref{NaSt_fullydisc}. We define interpolated velocity snapshots $\snapVelMod{j} = \lagrange{\hilbRefVel}\velh[j]\in\hilbRefVel$ and interpolated pressure snapshots $\snapPrsMod{j} = \lagrange{\hilbRefPrs}\prsh[j]\in\hilbRefPrs$ for $j=1,\dots,\nSnapshots$, where $\lagrange{\hilbRefVel}$ and $\lagrange{\hilbRefPrs}$ are Lagrange interpolation operators onto the velocity reference space $\hilbRefVel$ and the pressure reference space $\hilbRefPrs$, respectively. 

In order to define a POD of the interpolated velocity snapshots in the sense of \cref{sec:pod}, we specify a Hilbert space $\hilbSnap=\hilbVel$, a set of snapshot weights $\weight{1}=\dots=\weight{\nSnapshots}=\Delta t$ and a POD dimension $\nPodVel\leq\dim(\spn(\snapVelMod{1},\dots,\snapVelMod{\nSnapshots}))$. As a result, we obtain velocity POD basis functions $\podVel{1}\dots,\podVel{\nPodVel}\in\spn(\snapVelMod{1},\dots,\snapVelMod{\nSnapshots})$. For the pressure, we choose the space $\hilbSnap=\hilbPrs$, a set of weights $\weight{1}=\dots=\weight{\nSnapshots}=\Delta t$ and a POD dimension $\nPodPrs\leq\dim(\spn(\snapPrsMod{1},\dots,\snapPrsMod{\nSnapshots}))$. A POD provides pressure POD basis functions $\podPrs{1}\dots,\podPrs{\nPodPrs}\in\spn(\snapPrsMod{1},\dots,\snapPrsMod{\nSnapshots})$. The POD bases approximate the respective interpolated snapshots optimally in the sense of \eqref{minPOD}.

\subsection{Stabilization with supremizer functions}\label{sec:stabilization_supremizer}

The stability of the velocity-pressure reduced-order model provided by \eqref{tmp} depends on the choice of the subspaces $\hilbPodVel$ and $\hilbPodPrs$. We use $\hilbPodPrs=\spn(\podPrs{1}\dots,\podPrs{\nPodPrs})$ as a pressure subspace, with basis functions according to \cref{sec:supremizersPod}. In order to derive a stable reduced-order model, we define the velocity subspace
\begin{equation}\label{eq:hilbPodSup}
  \hilbPodSupVel[\nPod]:=\spn(\podVel{1},\dots,\podVel{\nPodVel},\podSup{1},\dots,\podSup{\nPodPrs}),
\end{equation}
where the basis functions $\podVel{1},\dots,\podVel{\nPodVel}$ are defined in \cref{sec:supremizersPod} and the additional basis functions $\podSup{1},\dots,\podSup{\nPodPrs}\in\hilbRefVel$ are chosen such that an inf-sup stability constraint can be verified.

We express the inf-sup stability constraint as follows: There exists a $\infsupRom>0$ defined by
\[
  \infsupRom := \infstack{\testPrs\in\hilbPodPrs}{\testPrs\neq 0}\supstack{\testVel\in\hilbPodSupVel[\nPod]}{\testVel\neq 0}\frac{\bilDivergence{\testVel}{\testPrs}}{\|\testVel\|_{\hilbVel}\|\testPrs\|_{\hilbPrs}}.
\]
In the following we show how suitable functions $\podSup{1},\dots,\podSup{\nPodPrs}$ can be found such that the stability condition is fulfilled \cite{BMQR15,RV07}.

We introduce a linear map $\supMap:\hilbPrs\rightarrow\hilbRefVel$ by the following problem: For given $\testPrs\in\hilbPrs$, find  $\supMap\testPrs\in\hilbRefVel$ such that
\begin{equation}\label{eq:StabilizerProblem}
  (\supMap\testPrs,\testVel)_{\hilbVel} = \bilDivergence{\testVel}{\testPrs}\quad\forall\testVel\in\hilbRefVel.
\end{equation}
From the Riesz representation theorem follows
\begin{equation}\label{eq:StabilizerMap}
  \|\supMap\testPrs\|_{\hilbRefVel} = \sup_{\substack{\testVel\in\hilbRefVel\\\testVel\neq 0}}\frac{\bilDivergence{\testVel}{\testPrs}}{\|\testVel\|_{\hilbVel}}
  = \sup_{\substack{\testVel\in\hilbRefVel\\\testVel\neq 0}}\frac{(\supMap\testPrs,\testVel)_{\hilbVel}}{\|\testVel\|_{\hilbVel}}\quad\forall\testPrs\in\hilbPrs,
\end{equation}
which means that
\begin{equation}\label{eq:ArgSup}
  \supMap\testPrs = \arg\sup_{\substack{\testVel\in\hilbRefVel\\\testVel\neq 0}}\frac{(\supMap\testPrs,\testVel)_{\hilbVel}}{\|\testVel\|_{\hilbVel}}\quad\forall\testPrs\in\hilbPrs.
\end{equation}
Therefore, we enrich the velocity space with functions $\podSup{i} = \supMap\podPrs{i}$ for $i=1,\dots,\nPodPrs$. 

To show that the stability condition is fulfilled for the proposed choice of enrichment functions, we relate the stability constant $\infsupRom$ of the reduced-order model with the stability constant $\infsupCom$ of the reference spaces $(\hilbRefVel,\hilbRefPrs)$, see \cite{Bab73,Bre74,Lad63}. We proceed like Proposition 2 of \cite{BMQR15} and Lemma 3.1 of \cite{RV07}:
\begin{align*}
  0
  <\infsupCom
  &:= \inf_{\substack{\testPrs\in\hilbRefPrs\\\testPrs\neq 0}}\;\sup_{\substack{\testVel\in\hilbRefVel\\\testVel\neq 0}}\frac{\bilDivergence{\testVel}{\testPrs}}{\|\testVel\|_{\hilbVel}\|\testPrs\|_{\hilbPrs}}
  \stackrel{\hilbPodPrs\subset\hilbRefPrs}{\leq}
  \infstack{\testPrs\in\hilbPodPrs}{\testPrs\neq 0}\;
  \supstack{\testVel\in\hilbRefVel}{\testVel\neq 0}\;
  \frac{\bilDivergence{\testVel}{\testPrs}}{\|\testVel\|_{\hilbVel}\|\testPrs\|_{\hilbPrs}}\\
  &\stackrel{\eqref{eq:StabilizerMap},\eqref{eq:ArgSup}}{=}\inf_{\substack{\testPrs\in\hilbPodPrs\\\testPrs\neq 0}}
  \frac{\bilDivergence{\supMap \testPrs}{\testPrs}}{\|\supMap \testPrs\|_{\hilbVel}\|\testPrs\|_{\hilbPrs}}\leq \infstack{\testPrs\in\hilbPodPrs}{\testPrs\neq 0}\;
  \supstack{\testVel\in\{\supMap \prs\colon \prs\in\hilbPodPrs\}}{\testVel\neq 0}\;
  \frac{\bilDivergence{w}{\testPrs}}{\|w\|_{\hilbVel}\|\testPrs\|_{\hilbPrs}}\\
  &\stackrel{\eqref{eq:hilbPodSup}}{\leq}
  \infstack{\testPrs\in\hilbPodPrs}{\testPrs\neq 0}\;
  \supstack{\testVel\in\hilbPodSupVel[\nPod]}{\testVel\neq 0}\;
  \frac{\bilDivergence{w}{\testPrs}}{\|w\|_{\hilbVel}\|\testPrs\|_{\hilbPrs}}=\infsupRom.
\end{align*}
The final relation $0 <\infsupCom\leq\infsupRom$ states that the POD model with an enriched velocity space is inf-sup stable as long as $(\hilbRefVel,\hilbRefPrs)$ is an inf-sup stable pair of spaces.

Note that \eqref{eq:StabilizerProblem} allows a computation of the stabilizer functions based on the pressure snapshots instead of the pressure POD. \Cref{sec:supremizersPod} states that $\podPrs{1}, \dots,$ $\podPrs{\nPodPrs}\in\spn(\snapPrsMod{1},\dots,\snapPrsMod{\nSnapshots})$. Consequently, for each $\iPod=1,\dots,\nPodPrs$ we can write
\begin{equation}\label{eq:StabilizerLinearCombination}
  \podSup{\iPod} = \supMap\podPrs{\iPod} = \supMap\sum_{j=1}^{\nSnapshots}\snapPrsMod{j}\xi_\iPod^j = \sum_{j=1}^{\nSnapshots}(\supMap\snapPrsMod{j})\xi_\iPod^j,
\end{equation}
where the coefficients $\xi_\iPod^1,\dots,\xi_\iPod^\nSnapshots$ can be obtained from the pressure POD computation. This means, in a first step, we compute $\supMap\snapPrsMod{1},\dots,\supMap\snapPrsMod{\nSnapshots}$ via \eqref{eq:StabilizerProblem}. This involves the solution of a linear system of equations on the reference finite element space for each pressure snapshot. In a second step, we compute $\podSup{1},\dots,\podSup{\nPodPrs}$ by linearly combining $\supMap\snapPrsMod{1},\dots,\supMap\snapPrsMod{\nSnapshots}$ according to \eqref{eq:StabilizerLinearCombination}. The result is equivalent to the supremizers obtained from the pressure POD basis functions.

\subsection{Complexity}

Concerning the complexity of constructing a reduced-order model, there are two main differences between the velocity ROM approach and the velocity-pressure ROM approach. First of all, in the velocity ROM approach a POD basis is computed only for the velocity variable, whereas in the velocity-pressure ROM approach an additional POD basis is computed for the pressure snapshots. Second, the construction of divergence-free POD modes in the velocity ROM approach requires the solution of \cref{problem:leraySaddlePoint}, whose complexity depends on $\text{dim}(\hilbRefVel)+\text{dim}(\hilbRefPrs)$. In contrast, the construction of the supremizer functions in the velocity-pressure ROM approach requires the solution of \cref{eq:StabilizerProblem}, whose complexity depends on $\text{dim}(\hilbRefVel)$. 

For the complexity of solving a velocity-pressure reduced-order model, we concentrate on the setup and solution of a linear system in a Newton step. The setup of the system matrix and right-hand side involves the third-order tensor originating from the convective terms. The solution  of the linear system requires the factorization of a dense matrix which contains a zero block in case of the velocity-pressure approach. As a result, the complexity of solving a velocity-pressure ROM is of order $\mathcal{O}((\nPodVel + \nPodPrs)^3)$ while the complexity of solving a velocity ROM is only of order $\mathcal{O}(\nPodVel^3)$.

\begin{remark}
 Since the computation of supremizer functions can be expensive in practical applications, we like to mention some variations and alternatives to the stabilization using supremizer enrichment. In \cite{BMQR15} an approximate supremizer computation is proposed which enables an efficient offline-online decomposition while preserving stability properties. An alternative to stabilization with supremizers is given in \cite{BBI09,CIJS14,WLBI10} as a velocity-pressure reduced-order model with a residual-based stabilization approach. A transfer of the supremizer stabilization technique to the context of finite volume approximations is carried out in \cite{SR18}, where a comparison to a stabilization based on a pressure Poisson equation in the online phase is provided.
\end{remark}

\section{Inhomogeneous Dirichlet data} \label{sec:inh_dirichlet}

\newcommand{\dirVel}{\check{\vel}}

\newcommand{\velMod}[1][]{\tilde y^{#1}}
\newcommand{\prsMod}[1][]{\tilde p^{#1}}
\newcommand{\dirMod}[1][]{\tilde\dir^{#1}}
\newcommand{\dirFeSpace}[1][]{V_\text{D}^{#1}}
\newcommand{\dirRef}{\tilde V_\text{D}}

So far, we have studied the incompressible Navier-Stokes problem with homogeneous Dirichlet boundary conditions. In the following, we extend the scope to problems involving inhomogeneous Dirichlet data. The main idea is to subtract a suitable chosen lifting function from the snapshot data leading to homogeneous Dirichlet boundary conditions for the reduced basis functions and then add the lifting function in the expansion of the velocity field. For PDEs with a single parametrized Dirichlet boundary this is referred to as \textit{control function method} in \cite{GPT99} and generalized to multiple parameters in \cite{GPS07,Ull14}.

In the context of POD-Galerkin modeling based on adaptive finite element snapshots, the main challenge is to find such lifting functions for each space-adapted snapshot. For the derivation of a velocity POD-Galerkin model, we further must ensure that these continuous extensions fulfill the correct weak divergence-free property.\\
Another alternative to handle inhomogeneous Dirichlet conditions is the \textit{penalty method} \cite{GPT99}, where the snapshot data is not homogenized but the inhomogeneous Dirichlet data is enforced in a weak form in the Galerkin projection. Moreover, in \cite{GPS07} a different approach is proposed which uses a modification of the POD basis utilizing a QR decomposition such that some of the reduced basis functions fulfill the homogeneous and some fulfill the inhomogeneous boundary data. A transfer of these approaches to the case of space-adapted snapshots is not carried out within the scope of this work.

We extend \eqref{NaSt} to the case of inhomogeneous Dirichlet boundary conditions by introducing Dirichlet boundary data $\velDir\colon[0,T)\times\partial\Omega\rightarrow \mathbb R$. The resulting problem reads as follows: Find a velocity field $\dirVel$ and a pressure field $\prs$ such that
\begin{subequations}\label{eq:dirichlet}
\begin{alignat}{4}
\dirVel_t + (\dirVel \cdot \nabla) \dirVel - Re^{-1} \Delta \dirVel + \nabla \prs &= f &\qquad& \text{in } (0,T) \times \Omega,\label{eq:dirichlet_momentum}\\
\divergence \; \dirVel &= 0 &\qquad& \text{in } (0,T) \times \Omega, \label{eq:dirichlet_continuity}\\
\dirVel &= \velDir      &\qquad& \text{in } (0,T)\times\partial\Omega,\label{eq:dirichlet_dirichlet}\\
\dirVel &= \velInit &\qquad& \text{in } \{0\}\times\Omega,\label{eq:dirichlet_initial}
\end{alignat}
\end{subequations}
where $\divergence\velInit = 0$ in $\Omega$ and $\velInit = \velDir(0)$ on $\partial\Omega$. In the following, we derive a homogenized version of this problem, which provides a foundation for the subsequent finite element discretization and reduced-order modeling.

\subsection{Homogenized equations}

We assume that the function $\velDir$ is sufficiently regular, so that it can be continuously extended by a function $\dir\colon[0,T)\times\bar\Omega\rightarrow \mathbb R$ with $\dir(t)|_{\partial\Omega}=\velDir(t)$ for $t\in(0,T)$. The regularity requirements on $\dir$ are such that a unique weak solution exists, see \cite{Raymond2007}. In \cref{sec:lifting} we provide a concrete choice of $\dir$ by computation.

We homogenize \eqref{eq:dirichlet} by subtracting the boundary function from the inhomogeneous velocity solution such that the homogeneous velocity field is given by $\vel = \dirVel - \dir$. Substituting $\dirVel$ in \eqref{eq:dirichlet}, we obtain the following homogenized problem: Find $\vel$ and $\prs$ such that
\begin{subequations}\label{eq:dirichlet_homogenized}
\begin{alignat}{4}
\notag\vel_t + (\vel \cdot \nabla) \vel + (\dir \cdot \nabla) \vel + (\vel \cdot \nabla) \dir  - Re^{-1} \Delta \vel + \nabla \prs \hspace{-5cm}\\ &= f - (\dir \cdot \nabla) \dir + Re^{-1} \Delta \dir - \dir_t && \qquad\text{in } (0,T) \times \Omega, \\
\divergence \; \vel &= -\divergence \; \dir  && \qquad \text{in } (0,T) \times \Omega, \\
\vel  &= 0      && \qquad \text{in } (0,T)\times\partial\Omega,\\
\vel  &= \velInit - \dir&& \qquad \text{in } \{0\}\times\Omega,
\end{alignat}
\end{subequations}
where $\divergence\velInit = 0$ in $\Omega$ and $\velInit = \dir(0)$ on $\partial\Omega$. We proceed like in \cref{sec:setting}, but now using the homogenized equations.

In order to derive a time-discrete weak form of the homogenized problem, we implement the time integrals involving the Dirichlet data using the right rule, which evaluates the Dirichlet data at the new time instance. For ease of notation, we define $f^j:=f(t^j)$ and $\dir[j] := \dir(\ti[j])$ for $j=0,\dots,\nSnapshots$. As a result, the time-discrete weak form of the homogenized problem consists in finding sequences $\vel[1],\dots,\vel[\nSnapshots]\in\hilbVel$ and $\prs[1],\dots,\prs[\nSnapshots]\in\hilbPrs$, for given $\vel[0]=\velInit - \dir[0]$ with $\velInit\in\Hdiv$, such that
\begin{subequations}\label{eq:dirichlet_timeDiscrete}
\begin{alignat}{2}
  \inner[\Big]{\frac{\vel[j]-\vel[j-1]}{\Delta t}}{v}{} + c(\vel[j],\vel[j],v) + c(\dir[j],\vel[j],v) + c(\vel[j],\dir[j],v)+ a(\vel[j],v) + b(v,\prs[j]) \hspace{-10cm}\notag\\
     & =   \duality{f^j}{v}{} - c(\dir[j],\dir[j],v) - a(\dir[j],v) - \inner[\Big]{\frac{\dir[j]-\dir[j-1]}{\Delta t}}{v}{}&\quad& \forall v\in \mathcal{V},\label{eq:dirichlet_timeDiscrete_momentum}\\
  b(\vel[j],q) & =  -b(\dir[j],q)&\quad& \forall q \in \mathcal{Q}\label{eq:dirichlet_timeDiscrete_continuity}
\end{alignat}
\end{subequations}
for $j=1,\dots,\nSnapshots$.

We utilize an adaptive finite element method, so that the fully discrete homogenized Navier-Stokes problem reads as follows: For given $\velh[0]=\velInit - \dir[0]$ with $\velInit\in\Hdiv$, find $\velh[1]\in\hilbSnapVel[1],\dots,\velh[\nSnapshots]\in\hilbSnapVel[\nSnapshots]$ and $\prsh[1] \in \hilbSnapPrs[1],\dots,\prsh[\nSnapshots] \in \hilbSnapPrs[\nSnapshots]$ such
that
\begin{subequations}\label{eq:dirichlet_fullyDiscrete}
\begin{alignat}{2}
  \inner[\Big]{\frac{\velh[j]-\velh[j-1]}{\Delta t}}{v}{} + c(\velh[j],\velh[j],v) + c(\dir[j],\velh[j],v) + c(\velh[j],\dir[j],v) + a(\velh[j],v) + b(v,\prsh[j]) \hspace{-10cm}\notag\\
    & = \duality{f^j}{v}{} - c(\dir[j],\dir[j],v) - a(\dir[j],v) - \inner[\Big]{\frac{\dir[j]-\dir[j-1]}{\Delta t}}{v}{}&\quad& \forall v\in \hilbSnapVel[j],\\
  b(\velh[j],q) & =  -b(\dir[j],q)&\quad& \forall q \in \hilbSnapPrs[j]\label{eq:dirichlet_fullyDiscrete_continuity}
\end{alignat}
\end{subequations}
for $j=1,\dots,\nSnapshots$.

\subsection{Lifting function} \label{sec:lifting}

Based on \eqref{eq:dirichlet_fullyDiscrete}, we compute approximations to the inhomogeneous solutions $\dirVel(\ti[j])$ of \eqref{eq:dirichlet} by adding the continuous extension of the Dirichlet data, i.e.\ $\dirVelh[j]:=\velh[j]+\dir[j]$ for $j=0,\dots,\nSnapshots$. Regardless of the choice of lifting functions $\dir[0],\dots,\dir[\nSnapshots]$, we can guarantee that $\dirVelh[0]$ fulfills the initial condition \eqref{eq:dirichlet_initial} and $\dirVelh[1],\dots,\dirVelh[\nSnapshots]$ fulfill the Dirichlet condition \eqref{eq:dirichlet_dirichlet} by construction. Nevertheless, in order to solve \eqref{eq:dirichlet_fullyDiscrete} numerically, concrete candidates of $\dir[0],\dots,\dir[\nSnapshots]$ must be fixed, at least implicitly. Our approach to reduced-order modeling is not restricted to a particular choice. In the following, we provide suitable candidates which can be realized without the need to modify usual finite element codes.

Note that for the velocity finite element spaces it holds $\hilbSnapVel[j] \subset \hilbVel = H_0^1(\Omega)$. Thus, for the context of inhomogeneous Dirichlet conditions, we start by introducing the spaces $\dirFeSpace[j]$ for $j=1, \dots, \nSnapshots$, which denote the spaces spanned by the union of the finite element basis functions of $\hilbSnapVel[j]$ and the finite element basis functions associated with the corresponding Dirichlet boundary nodes. We assume that in \eqref{eq:dirichlet_fullyDiscrete} the integrals involving $\dir[j]$ are approximated by a numerical quadrature which consists of substituting the Lagrange interpolation of $\dir[j]$ onto $\dirFeSpace[j]$ and integrating the resulting piecewise polynomials exactly. We assume that by a finite number of refinements of any $\dirFeSpace[1],\dots,\dirFeSpace[\nSnapshots]$ one can find a common reference finite element space $\dirRef$ such that $\dirFeSpace[1],\dots, \dirFeSpace[\nSnapshots] \subset \dirRef$. Now, for $j=1,\dots,\nSnapshots$, we define lifting functions $\dir[j]$ as a sufficiently smooth continuous extension of the Dirichlet data $\velDir(\ti[j])$ into the domain $\Omega$ such that $\dir[j]$ is zero at all nodes of the reference finite element space $\hilbRefVel$. This is equivalent to the standard approach of using an \emph{approximate} Dirichlet lifting given by a Lagrangian interpolation of the Dirichlet data onto the finite element space at the boundary and a subsequent continuous extension using the finite element space in the interior, because we have
\[
  \begin{cases}
    \dir[j] = \velDir(t)&\text{ at all Dirichlet nodes of } \dirFeSpace[j],\\
    \dir[j] = 0&\text{ at all interior nodes of } \dirFeSpace[j],
  \end{cases}
\]
for $j=1,\dots,n$. A disadvantage of the standard approach is that it implies a Dirichlet lifting which satisfies the boundary data only in an approximate sense. Our description, on the other hand, delivers an output which is exact at the boundary. In particular, we have
\[
  \begin{cases}
    \dirVelh[j] = \velh[j]&\text{ at all interior nodes of } \dirRef\\
    \dirVelh[j] = \dir[j]&\text{ at all Dirichlet nodes of } \dirRef\\
    \velh[j] = 0&\text{ at all Dirichlet nodes of } \dirRef.
  \end{cases}
\]
This holds for all $\dirRef$ which fulfill our assumptions, without the need to specify a concrete candidate of $\dirRef$ during the adaptive finite element simulation. When the adaptive finite element simulation is finished and $\dirFeSpace[1],\dots,\dirFeSpace[\nSnapshots]$ are available, some $\dirRef$ can be computed by refinement and $\dir[j]$ can be evaluated at all nodal points of $\dirRef$. Therefore, we are even able to formulate a finite element discretization of \eqref{eq:dirichlet_timeDiscrete} on $(\hilbRefVel,\hilbRefPrs)$ using the same $\dir[0],\dots,\dir[\nSnapshots]$ as in \eqref{eq:dirichlet_fullyDiscrete}. Moreover, we are able to solve \eqref{eq:dirichlet_timeDiscrete} on subspaces of $(\hilbRefVel,\hilbRefPrs)$ using the same $\dir[0],\dots,\dir[\nSnapshots]$ as in \eqref{eq:dirichlet_fullyDiscrete}.

\begin{remark}
  In principle, it is possible to impose a weak divergence-free constraint on the homogenized velocity finite element solution by using lifting functions which are computed such that $b(\dir[j],q) = 0$ for all $q \in \hilbSnapPrs[j]$ and $j=0,\dots,\nSnapshots$. But this would require the solution of an additional stationary finite element problem for each $\dir[j]$. Moreover, this would not automatically imply a weak divergence-free property with respect to a reference pressure space $\hilbRefPrs$. An alternative to the implicit choice of the Dirichlet lifting function is its explicit choice at the level of the strong formulation \eqref{eq:dirichlet_homogenized}. Disadvantages would be a possibly larger support of such a lifting function and the effort of actually finding a suitable function. Also in this case, it would be attractive to impose a strong divergence-free constraint on $\dir[j]$, because this implies a weakly divergence-free homogenized velocity field $\velh[j]$. Still, finding a suitable candidate may be challenging in general.
\end{remark}

\subsection{Velocity POD-Galerkin model}\label{subsection:Velocity_pod_inh}

In the following, we derive a reduced-order model for the velocity field, based on the semi-discretized problem \eqref{eq:dirichlet_timeDiscrete}. We introduce the projections $\lerayInh[{\dir}]^{\hilbVel,\hilbPrs}$ according to \cref{problem:lerayMinimization} by
\begin{problem} \label{problem:lerayMinimizationDir}
For given $\snap \in \hilbSnap$, sufficiently smooth $\dir$ and given spaces $\hilbVel$ and $\hilbPrs$, find $\lerayInh[{\dir}]^{\hilbVel,\hilbPrs}(\snap)=\lerSnap\in\hilbVel$ which solves
\[
  \min_{v\in\hilbVel}\frac12\|v-\snap\|_{\hilbSnap}^2\quad\text{subject to}\quad b(v+\dir,\testPrs)=0\quad\;\forall \testPrs\in\hilbPrs.
\]
\end{problem}
We use $\hilbSnap=\hilbVel$ in the following. By definition of this projection onto a divergence-free space, we have $b(\lerayInh[{\dir[j]}]^{\hilbVel,\hilbPrs}(0)+\dir[j],\testPrs) = 0$ for all $\testPrs\in \hilbPrs$ and $j=0,\dots,\nSnapshots$. In \eqref{eq:dirichlet_timeDiscrete}, we substitute $\dir[j]=\dir[j]+\lerayInh[{\dir[j]}]^{\hilbVel,\hilbPrs}(0)-\lerayInh[{\dir[j]}]^{\hilbVel,\hilbPrs}(0)$ and reformulate the equations so that we can set $\velMod[j] = \vel[j] - \lerayInh[{\dir[j]}]^{\hilbVel,\hilbPrs}(0)$. We obtain the following problem, which is equivalent to \eqref{eq:dirichlet_timeDiscrete}: For given $\velMod[0]=\velInit - \dirMod[0]$ with $\velInit\in\Hdiv$ and $\dirMod[j] = \dir[j] + \lerayInh[{\dir[j]}]^{\hilbVel,\hilbPrs}(0)$ for $j=0,\dots,\nSnapshots$, find $\velMod[1],\dots,\velMod[\nSnapshots]\in\hilbVel$ and $\prs[1],\dots,\prs[\nSnapshots]\in\hilbPrs$ such that
\begin{subequations}\label{eq:dirichlet_Modified}
\begin{alignat}{4}
  \inner[\Big]{\frac{\velMod[j]-\velMod[j-1]}{\Delta t}}{v}{} + c(\velMod[j],\velMod[j],v) + c(\dirMod[j],\velMod[j],v) + c(\velMod[j],\dirMod[j],v) + a(\velMod[j],v) + b(v,\prs[j]) \hspace{-10cm}\notag\\
     & =  \duality{f^j}{v}{} - c(\dirMod[j],\dirMod[j],v) - a(\dirMod[j],v) - \inner[\Big]{\frac{\dirMod[j]-\dirMod[j-1]}{\Delta t}}{v}{}&\quad& \forall v\in \hilbVel,\\
  b(\velMod[j],q) & =  0&\quad& \forall q \in \hilbPrs
\end{alignat}
\end{subequations}
for $j=1,\dots,\nSnapshots$.

It can be shown that \eqref{eq:dirichlet_fullyDiscrete} is a discretization of \eqref{eq:dirichlet_Modified} by replacing $(\hilbVel,\hilbPrs)$ with $(\hilbSnapVel[j],\hilbSnapPrs[j])$ for $j=1,\dots,\nSnapshots$ in \eqref{eq:dirichlet_Modified}. Then, for the resulting solution holds $\velMod[j]=\velh[j] - \lerayInh[{\dir[j]}]^{\hilbSnapVel[j],\hilbSnapPrs[j]}(0)$ if $\dirMod[j] = \dir[j] + \lerayInh[{\dir[j]}]^{\hilbSnapVel[j],\hilbSnapPrs[j]}(0)$ for $j=0,\dots,\nSnapshots$. In this way, \eqref{eq:dirichlet_fullyDiscrete} can be used to obtain approximate solutions of \eqref{eq:dirichlet_Modified}.

We base the model equation on a discretization of \eqref{eq:dirichlet_Modified} using the pair of reference spaces $(\hilbRefVel,\hilbRefPrs)$ as test spaces together with $\dirMod[j] = \dir[j] + \lerayInh[{\dir[j]}]^{\hilbRefVel,\hilbRefPrs}(0)$. The resulting solutions are approximations to the solutions of \eqref{eq:dirichlet_Modified} using the original pair of spaces $(\hilbVel,\hilbPrs)$. We have shown that the solutions of \eqref{eq:dirichlet_Modified} using the original pair of spaces $(\hilbVel,\hilbPrs)$ are approximated by $\velh[j] - \lerayInh[{\dir[j]}]^{\hilbSnapVel[j],\hilbSnapPrs[j]}(0)$ for $j=1,\dots,\nSnapshots$ resulting from \eqref{eq:dirichlet_fullyDiscrete}. But these solutions are not weakly divergence-free with respect to the reference pair of spaces $(\hilbRefVel,\hilbRefPrs)$. Therefore, we have to modify them.

Following the velocity-ROM approach based on projected snapshots in \cref{sec:approach_i} we substitute $\velh[j] - \lerayInh[{\dir[j]}]^{\hilbSnapVel[j],\hilbSnapPrs[j]}(0)$ by their approximations $\lerayInh[{\dir[j]}]^{\hilbRefVel,\hilbRefPrs}(\velh[j]) - \lerayInh[{\dir[j]}]^{\hilbRefVel,\hilbRefPrs}(0)$ for $j=1,\dots,\nSnapshots$. Using now $\lerayInh[{\dir[j]}]^{\hilbRefVel,\hilbRefPrs}(\velh[j]) - \lerayInh[{\dir[j]}]^{\hilbRefVel,\hilbRefPrs}(0)$ as snapshots in a POD yields POD basis functions
\[
  \podVel{i}\in\spn\big(\lerayInh[{\dir[1]}](\velh[1]) - \lerayInh[{\dir[1]}](0),\dots,\lerayInh[{\dir[\nSnapshots]}](\velh[\nSnapshots]) - \lerayInh[{\dir[\nSnapshots]}](0)\big)\subset \hilbRefVelDiv\qquad\forall i=1,\dots,\nPodVel
\]
for some $\nPodVel\leq\nSnapshots$, which define a POD space $\hilbPodVel := \spn(\podVel{1}, \dots, \podVel{\nPodVel})\subset\hilbRefVelDiv$. 

In the time-discrete equation \eqref{eq:dirichlet_Modified}, we use the pair $(\hilbPodVel, \hilbRefPrs)$ as test and trial spaces. Consequently, the continuity equation is fulfilled by construction. For the pressure term, we have $b(v,\prs[j]) = 0$ for all $v\in\hilbPodVel$ and all $\prs[j]\in\hilbRefPrs$. The resulting reduced-order model is given by the following set of equations: For $\velPOD[0] = \velInit - \dirMod[0]$ with $\velInit\in\Hdiv$ and $\dirMod[j] = \dir[j] + \lerayInh[{\dir[j]}]^{\hilbRefVel,\hilbRefPrs}(0)$ for $j=0,\dots,\nSnapshots$, find $\velPOD[1],\dots,\velPOD[\nSnapshots]\in\hilbPodVel$ such that
\begin{align}\label{eq:dirichlet_velRom}
 & \inner[\Big]{\frac{\velPOD[j]-\velPOD[j-1]}{\Delta t}}{v}{} + c(\velPOD[j],\velPOD[j],v) + c(\dirMod[j],\velPOD[j],v) + c(\velPOD[j],\dirMod[j],v)     +a(\velPOD[j],v) \nonumber \\
 & \quad =  \duality{f^j}{v}{} -  c(\dirMod[j],\dirMod[j],v) - a(\dirMod[j],v) -  \inner[\Big]{\frac{\dirMod[j]-\dirMod[j-1]}{\Delta t}}{v}{} \quad \forall v\in \hilbPodVel
\end{align}
for $j=1,\dots,n$.
\begin{remark}
 Concerning the computational complexity, we have to additionally consider the projections $\lerayInh[{\dir[j]}]^{\hilbRefVel,\hilbRefPrs}(0)$ for $j=0, \dots, \nSnapshots$ in comparison to the homogeneous case. Therefore, the solution of \cref{problem:lerayMinimization} has to be computed $\nSnapshots+1$ times additionally to the projections of the homogeneous solutions $\velh[j]$.
\end{remark}
Following the velocity-ROM approach based on projected POD basis functions in \cref{sec:approach_ii}, we first introduce a set of modified homogeneous solutions $\snapVelMod{j} - \lerayInh[{\dir[j]}]^{\hilbRefVel,\hilbRefPrs}(0)$, for $j=1,\dots,\nSnapshots$. These modified snapshots can be constructed using e.g. a Lagrange interpolation of the original homogeneous solutions $\velh[j]$ onto the reference space $\hilbRefVel$ and an approximation of $\lerayInh[{\dir[j]}]^{\hilbSnapVel[j],\hilbSnapPrs[j]}(0)$ for $j=1, \dots, \nSnapshots$. From these modified snapshots, we compute a POD basis $\podVelMod{1}, \dots, \podVelMod{\nPodVel} \in \hilbRefVel$. Note that these modes are in general not divergence-free. Thus, they are then projected onto the space $\hilbRefVelDiv$ by solving \cref{problem:leraySaddlePoint} with $\dir=0$. This leads to a divergence-free velocity POD space $\hilbPodVel = \spn\{ \podVel{1}, \dots, \podVel{\nPodVel} \} \subset \hilbRefVelDiv$. Replacing $(\hilbVel,\hilbPrs)$ by the pair $(\hilbPodVel,\hilbRefPrs)$ in \eqref{eq:dirichlet_Modified} leads to a reduced-order model of the form \eqref{eq:dirichlet_velRom}.

\subsection{Velocity-pressure POD-Galerkin model} \label{subsection:Velocity_pressure_pod_inh}

To derive a velocity-pressure reduced-order model of the homogenized problem \eqref{eq:dirichlet_fullyDiscrete}, we require a suitable inf-sup stable pair of reduced spaces. Since the homogenization does not alter the bilinear form $b(\cdot,\cdot)$, the inf-sup stability criterion stays the same. Therefore, we compute a pressure reduced space $\hilbPodPrs$ and a velocity reduced space $\hilbPodVel$ like in \cref{sec:supremizers}, but using Lagrange-interpolated velocity and pressure snapshots of \eqref{eq:dirichlet_fullyDiscrete} instead of \eqref{NaSt_fullydisc}. We derive a stable POD-Galerkin model from the time-discrete problem \eqref{eq:dirichlet_timeDiscrete} by using the pair $(\hilbPodPrs,\hilbPodVel)$ as test and trial spaces.

We solve the reduced-order model for the POD approximations $\velPOD[1],\dots,\velPOD[\nSnapshots]$ of the homogeneous velocity fields and the POD approximations $\prsPOD[1],\dots,\prsPOD[\nSnapshots]$ of the pressure fields. Finally, $\velPOD[j]+\dir[j]$ is a time-discrete reduced-order approximation of the velocity solution of the inhomogeneous problem.

\section{Numerical Example} \label{sec:example}

We use a regularized lid-driven flow in a cavity as a numerical example. It describes the evolution of a flow in a confined domain $\Omega=(0,1)\times(0,1)$ with boundary $\partial\Omega$ during the time interval $[0,1]$. The governing equations are the incompressible Navier-Stokes equations with inhomogeneous Dirichlet data, as provided by \eqref{eq:dirichlet}. We specify a Reynolds number $Re=100$ and set $f(t,x)=0$. We choose a regularized lid-velocity according to \cite{DefrutosEA2016,John2016}, which has a quadratic velocity profile in $x_1$-direction and leads to a smooth transition to a zero velocity at the upper corners. In particular, it implies $\divergence \vel = 0$ in the corners and avoids a singularity of the pressure in the upper right corner. An additional regularization in time allows for a smooth startup from $\velInit(x)=0$. Hence, the Dirichlet data is given by $\velDir(t,x) = \velDir^t(t)\velDir^x(x)$ for all $(t,x)\in[0,1]\times\partial\Omega$, where
\begin{align*}
  \velDir^t(t) &=
    \begin{cases}
      1-\frac14(1-\cos((0.1-t)\pi/0.1))^2&\text{if }t\in[0,0.1),\\
      1 &\text{if }t\in[0.1,1],
    \end{cases}\\
  \velDir^x(x) &=
    \begin{cases}
      1-\frac14(1-\cos((0.1-x_1)\pi/0.1))^2&\text{if }x_2=1, x_1\in[0,0.1],\\
      1&\text{if }x_2=1, x_1\in(0.1,0.9), \\
      1-\frac14(1-\cos((x_1-0.9)\pi/0.1))^2&\text{if }x_2=1, x_1\in[0.9,1],\\
      0&\text{otherwise}.
  \end{cases}
\end{align*}

\subsection{Discretization}

We discretize the example problem using a space-adaptive extension of our Matlab finite element code \cite{Ull16}. The initial finite element mesh $\initialtriangulation$ is given by a criss-cross triangulation of a $8\times 8$ square pattern, see \cref{fig:mesh} on the left. We choose $\varepsilon = 0.01$ as a stopping tolerance and $\theta = 0.1$ as a refinement parameter in the adaptive \cref{Alg:AdFE}. For the time discretization, we set the number of discrete time intervals equal to $\nSnapshots=100$, so that $\Delta t = 0.01$.

We run the fully discrete Navier-Stokes problem \eqref{eq:dirichlet_fullyDiscrete} with the provided discretization parameters to compute a set of velocity and pressure solutions. \Cref{fig:solution} presents the components of the adaptive finite element solution at times $\ti = 0.1, 0.3, 1.0$ as well as the corresponding adapted finite element meshes.

\newenvironment{MyPicture}[1][]{
  \begin{axis}[
    title style={yshift=-2mm},
    width = 31.45mm,
    height = 31.45mm,
    point meta min = \minVelY,
    point meta max = \maxVelY,
    axis x line=none,
    axis y line=none,
    enlargelimits=false,
    scale only axis,
    #1]}{
  \end{axis}}

\newenvironment{MyPictureWithoutColorbar}[1][]{
  \begin{pgfinterruptboundingbox}
  \begin{MyPicture}[#1]}{
  \end{MyPicture}
  \end{pgfinterruptboundingbox}
  \useasboundingbox(-0cm,-0.2125cm)rectangle(3.4cm,3.825cm);}

\newenvironment{MyPictureWithColorbar}[1][]{
  \begin{pgfinterruptboundingbox}
  \begin{MyPicture}[#1]}{
  \end{MyPicture}
  \end{pgfinterruptboundingbox}
  \useasboundingbox(-0cm,-0.2125cm)rectangle(5.1cm,3.825cm);}

\newcommand{\minVelY}{0}
\newcommand{\maxVelY}{0}

\begin{figure}
  \begin{center}
    \ifUseExternalPngs
      \includegraphics{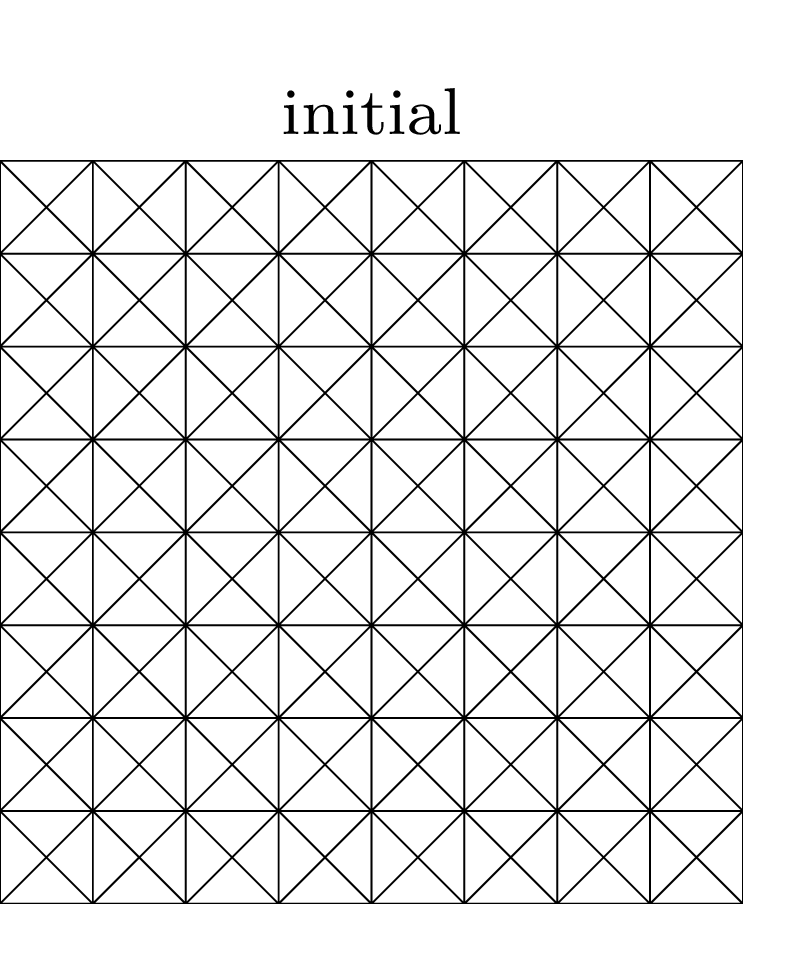}
    \else
      \begin{tikzpicture}
        \begin{MyPictureWithoutColorbar}[title = initial]
          \input{main15_mesh_1.tex}.
        \end{MyPictureWithoutColorbar}
      \end{tikzpicture}
    \fi
    \hspace{2mm}
    \ifUseExternalPngs
      \includegraphics{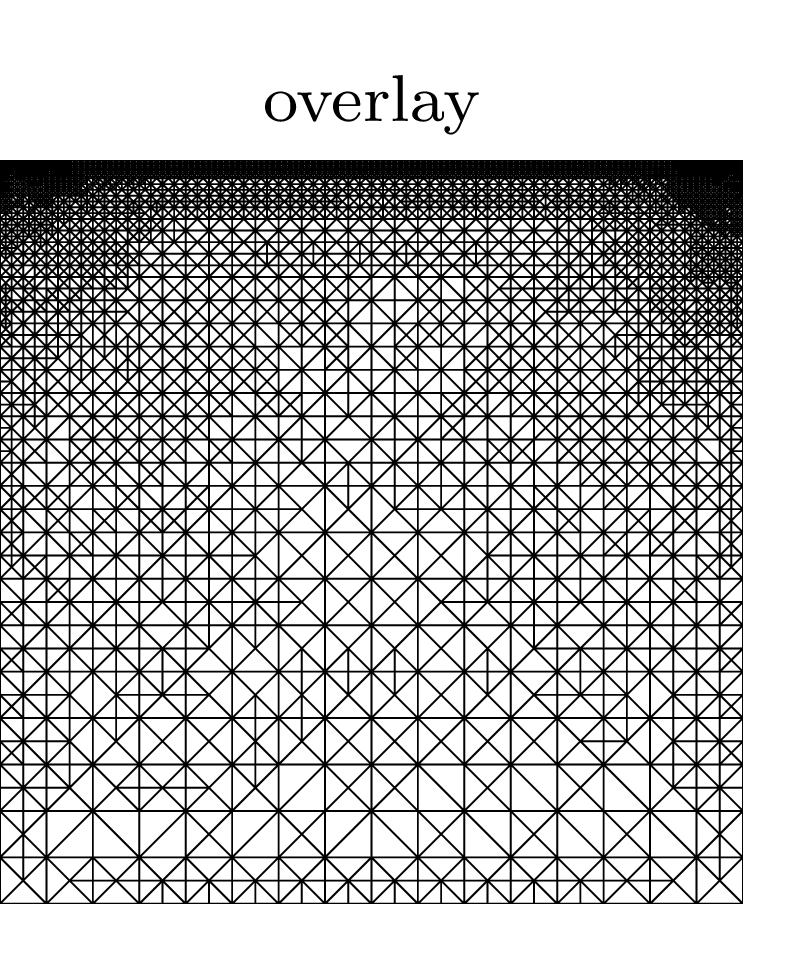}
    \else
      \begin{tikzpicture}
        \begin{MyPictureWithoutColorbar}[title = overlay]
          \input{main15_mesh_common.tex}
        \end{MyPictureWithoutColorbar}
      \end{tikzpicture}
    \fi
  \end{center}
  \caption{Initial triangulation $\initialtriangulation$ (left) and overlay of all adapted triangulations (right).\label{fig:mesh}}
\end{figure}

\begin{figure}
\begin{center}
  \ifUseExternalPngs
    \includegraphics{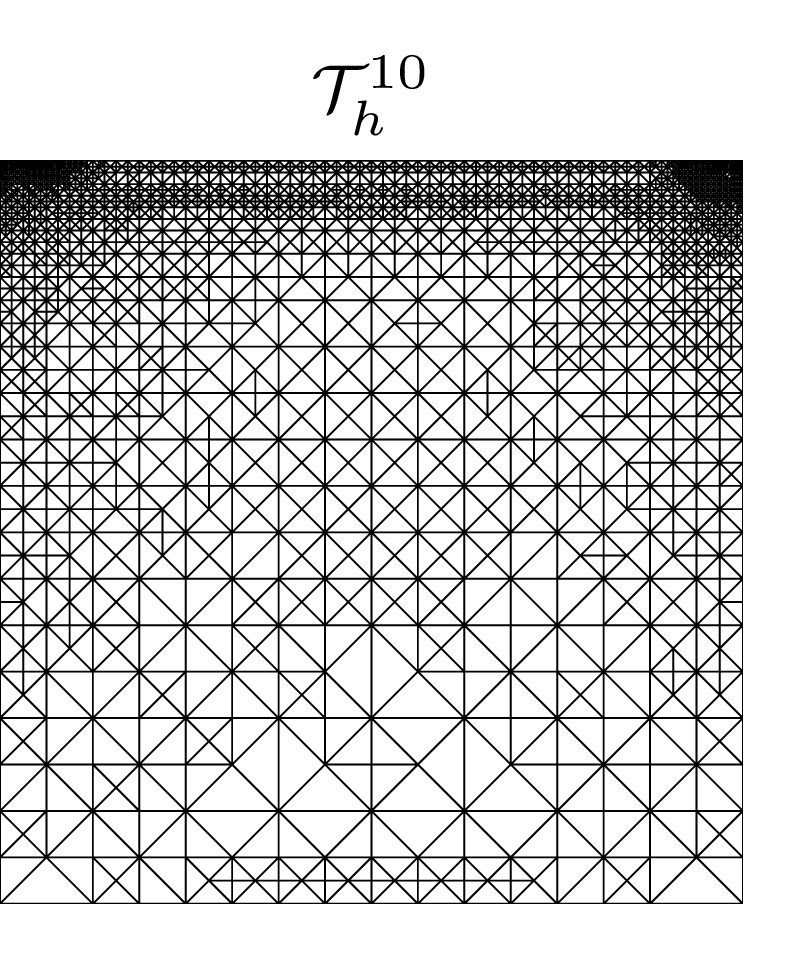}
  \else
    \begin{tikzpicture}
      \begin{MyPictureWithoutColorbar}[title = {$\triangulation{10}$}]
        \input{main15_mesh_11.tex}
      \end{MyPictureWithoutColorbar}
    \end{tikzpicture}
  \fi
\ifUseExternalPngs
    \includegraphics{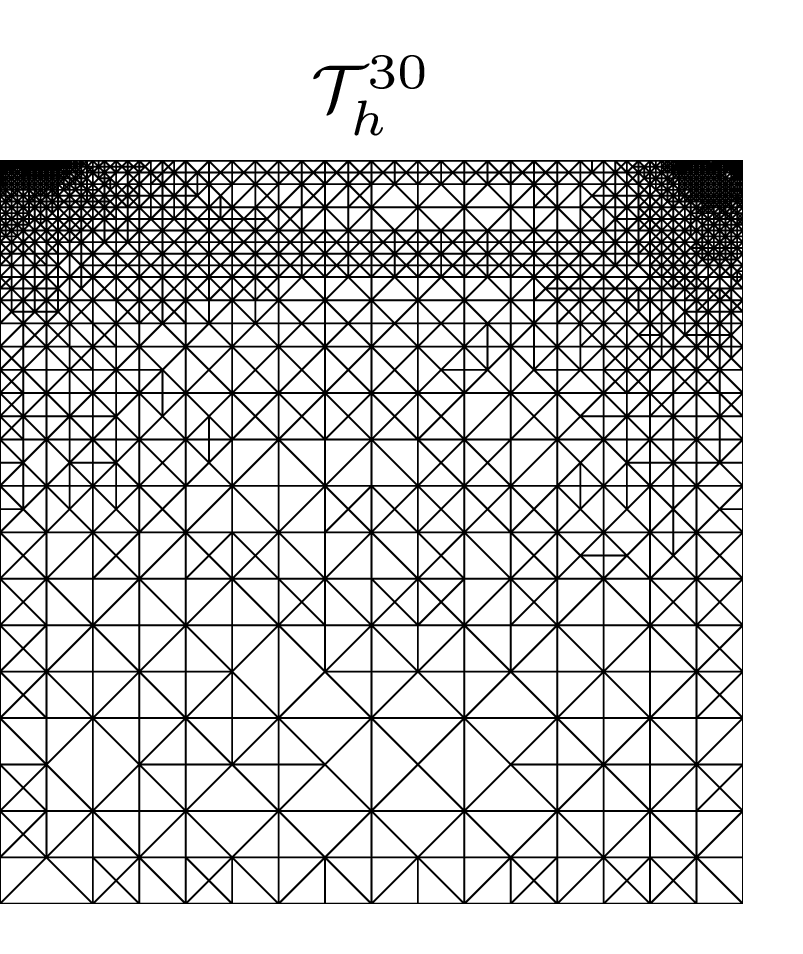}
  \else
    \begin{tikzpicture}
      \begin{MyPictureWithoutColorbar}[title = {$\triangulation{30}$}]
        \input{main15_mesh_31.tex}
      \end{MyPictureWithoutColorbar}
    \end{tikzpicture}
  \fi
 \ifUseExternalPngs
    \includegraphics{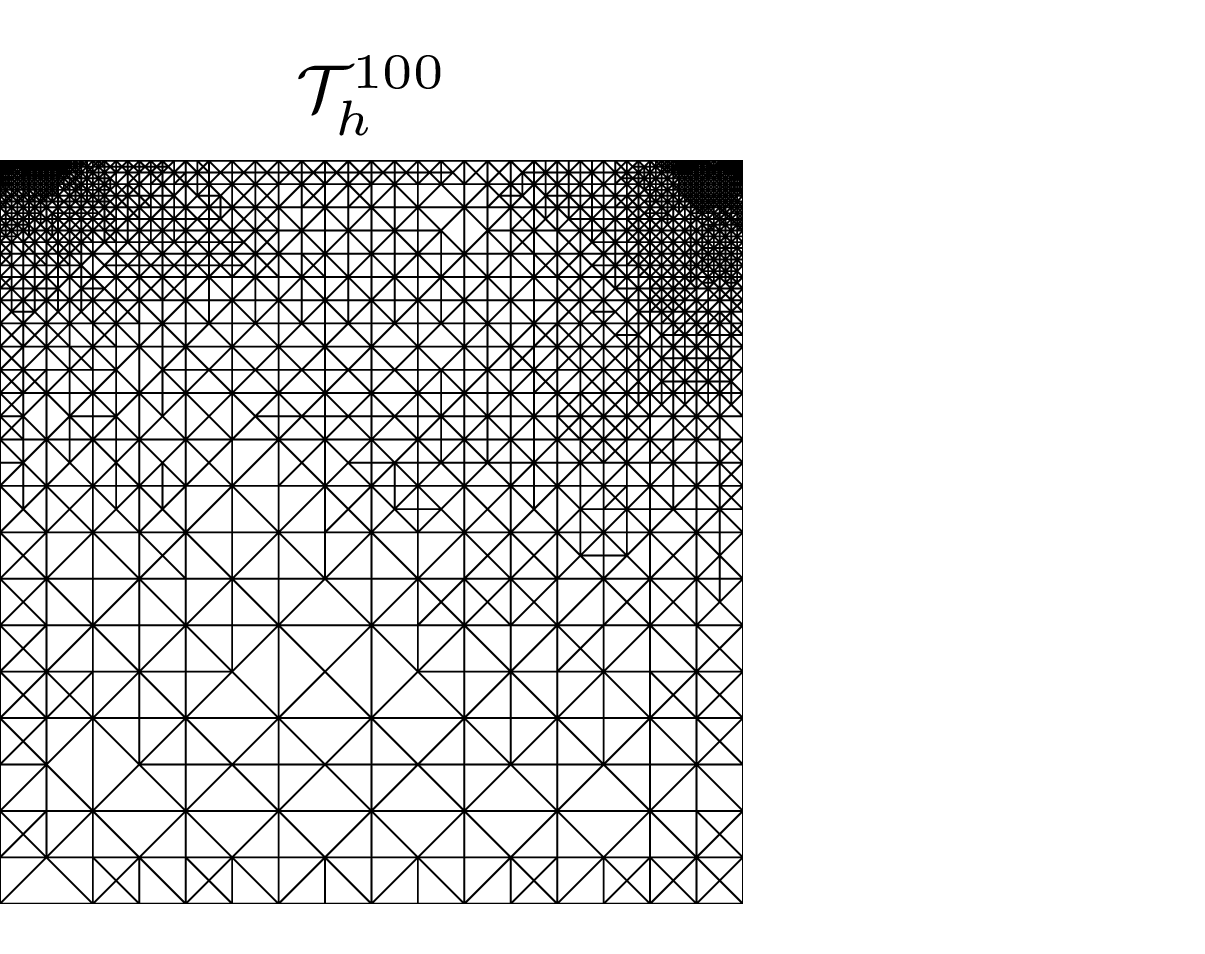}
  \else
    \begin{tikzpicture}
      \begin{MyPictureWithColorbar}[title = {$\triangulation{100}$}]
        \input{main15_mesh_101.tex}
      \end{MyPictureWithColorbar}
    \end{tikzpicture}
  \fi
\\
  \renewcommand{\minVelY}{-1.0}
  \renewcommand{\maxVelY}{1.0}
  \ifUseExternalPngs
    \includegraphics{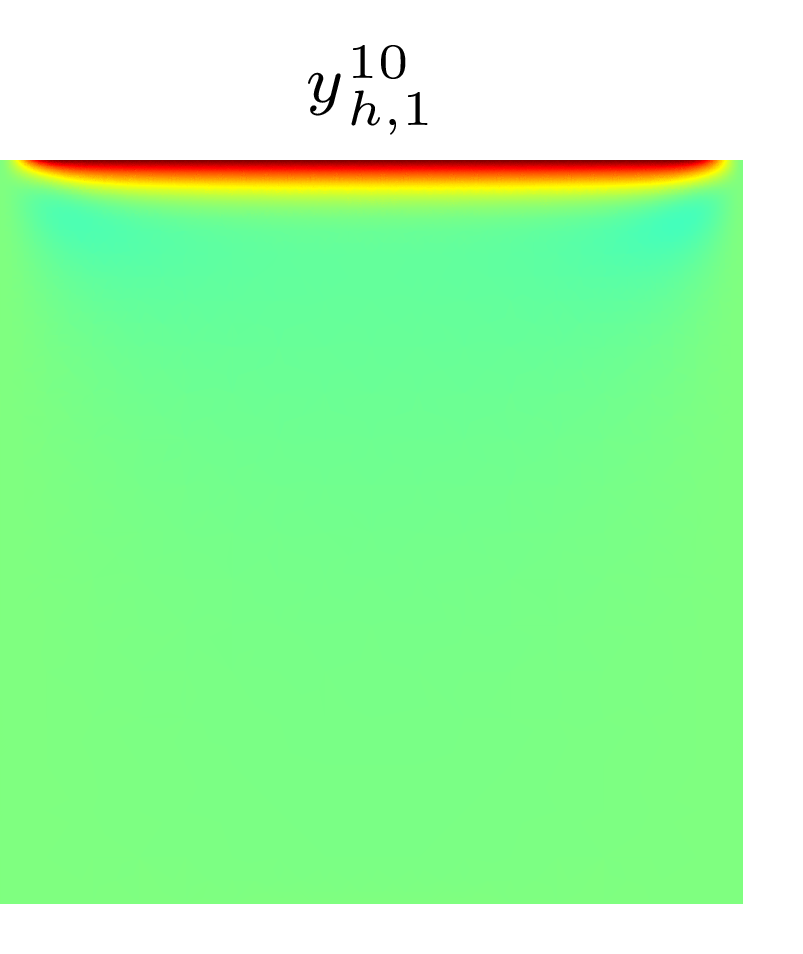}
  \else
  \begin{tikzpicture}
    \begin{MyPictureWithoutColorbar}[title = {$\velhx[10]$}]
      \input{main15_u_11.tex}
    \end{MyPictureWithoutColorbar}
  \end{tikzpicture}
  \fi
  \ifUseExternalPngs
    \includegraphics{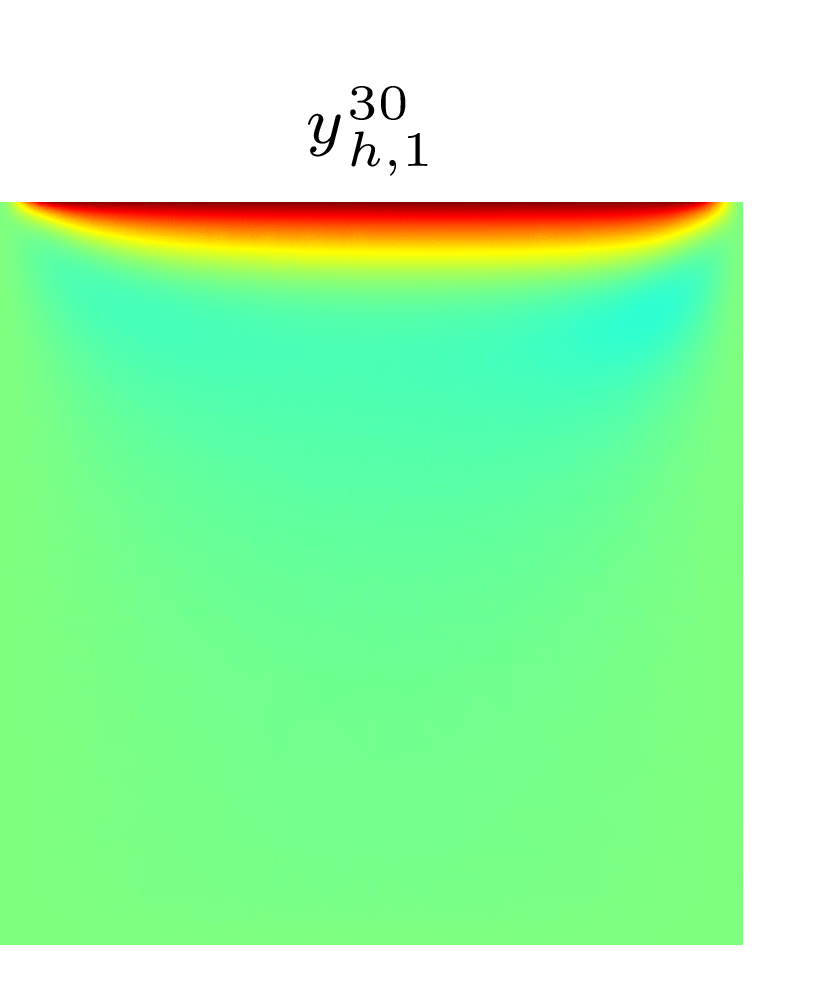}
  \else
    \begin{tikzpicture}
      \begin{MyPictureWithoutColorbar}[title = {$\velhx[30]$}]
        \input{main15_u_31.tex}
      \end{MyPictureWithoutColorbar}
      \useasboundingbox(-0cm,-0.25cm)rectangle(3.5cm,4cm);
    \end{tikzpicture}
  \fi
  \ifUseExternalPngs
    \includegraphics{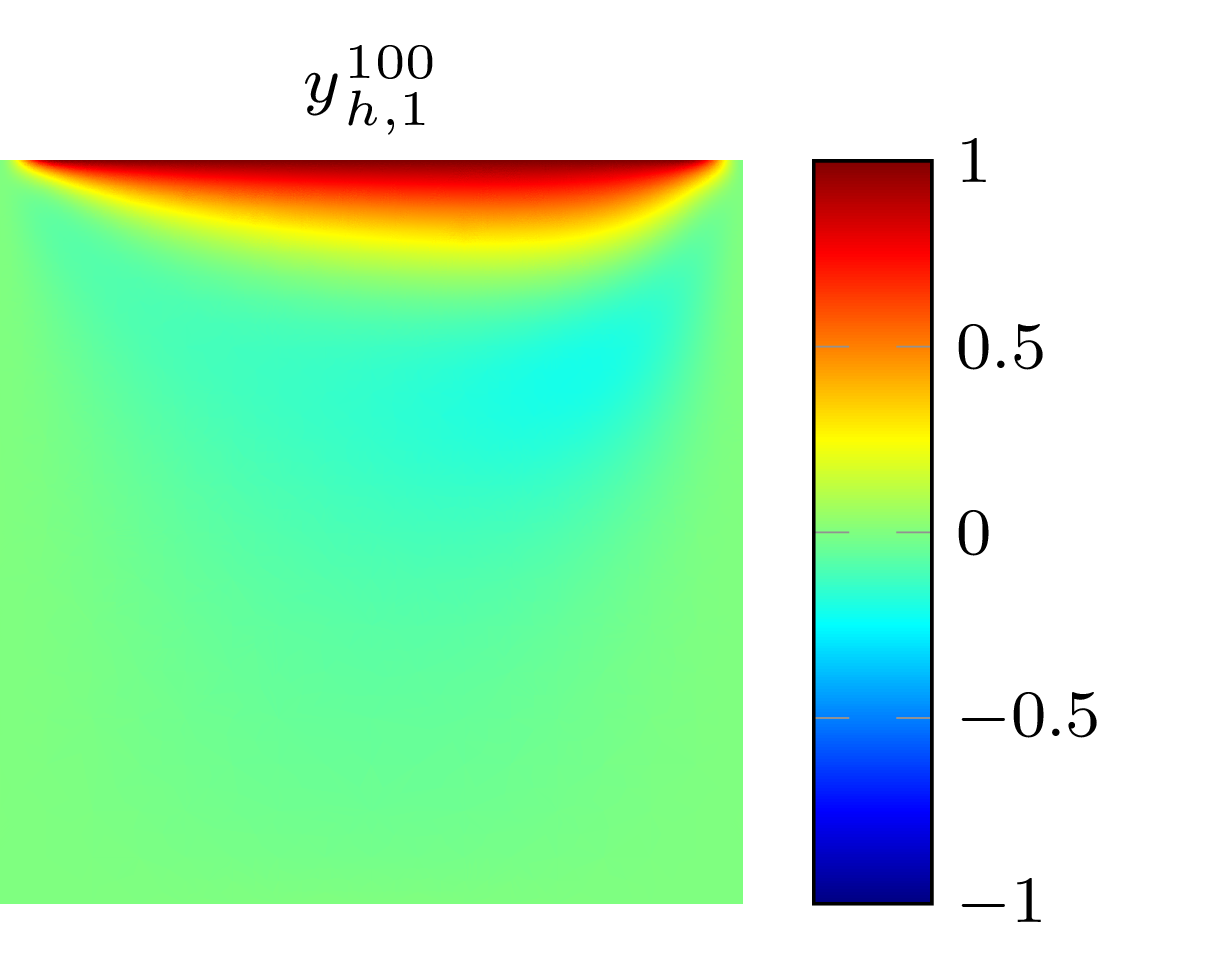}
  \else
    \begin{tikzpicture}
      \begin{MyPictureWithColorbar}[title = {$\velhx[100]$},colorbar,colorbar style = {ytick={-1,-0.5,0,0.5,1}}]
        \input{main15_u_101.tex}
      \end{MyPictureWithColorbar}
    \end{tikzpicture}
  \fi
\\
  \renewcommand{\minVelY}{-0.5}
  \renewcommand{\maxVelY}{0.5}
  \ifUseExternalPngs
    \includegraphics{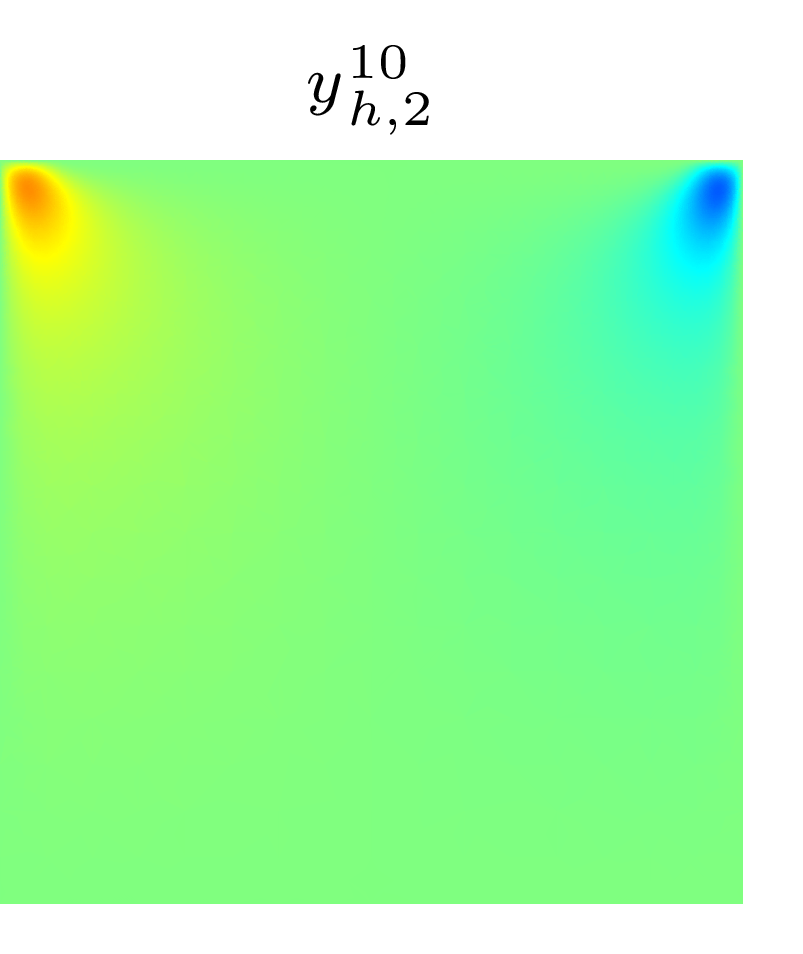}
  \else
    \begin{tikzpicture}[baseline, trim axis left]
      \begin{MyPictureWithoutColorbar}[title = {$\velhy[10]$}]
        \input{main15_v_11.tex}
      \end{MyPictureWithoutColorbar}
    \end{tikzpicture}
  \fi
  \ifUseExternalPngs
    \includegraphics{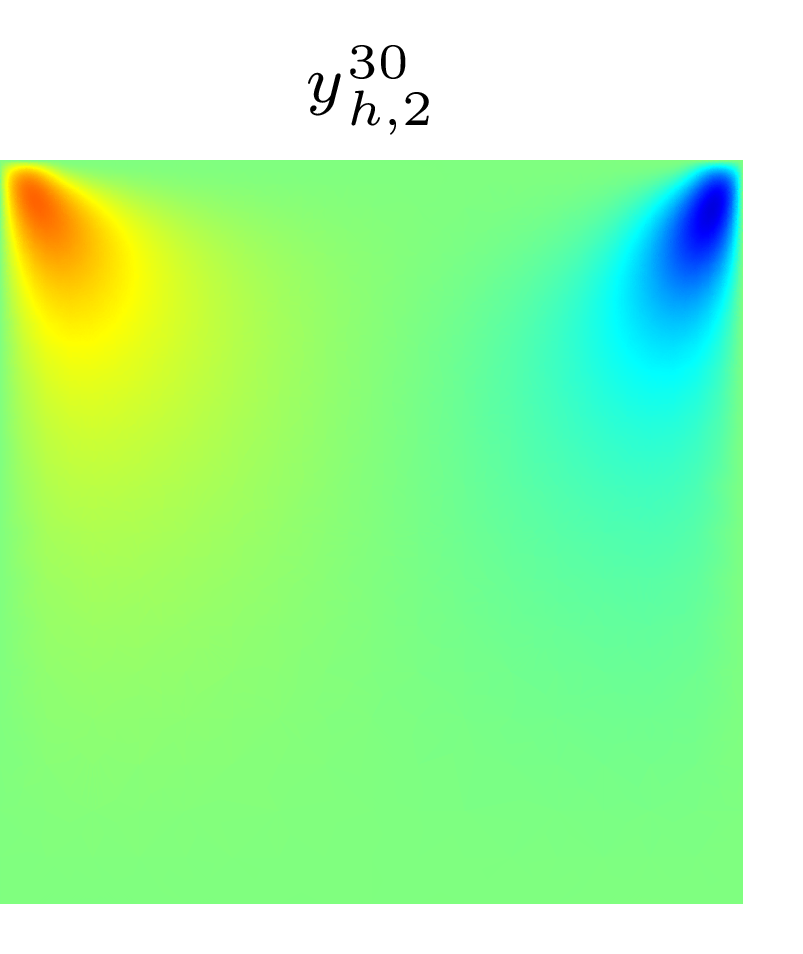}
  \else
    \begin{tikzpicture}[baseline, trim axis left]
      \begin{MyPictureWithoutColorbar}[title = {$\velhy[30]$}]
        \input{main15_v_31.tex}
      \end{MyPictureWithoutColorbar}
    \end{tikzpicture}
  \fi
  \ifUseExternalPngs
    \includegraphics{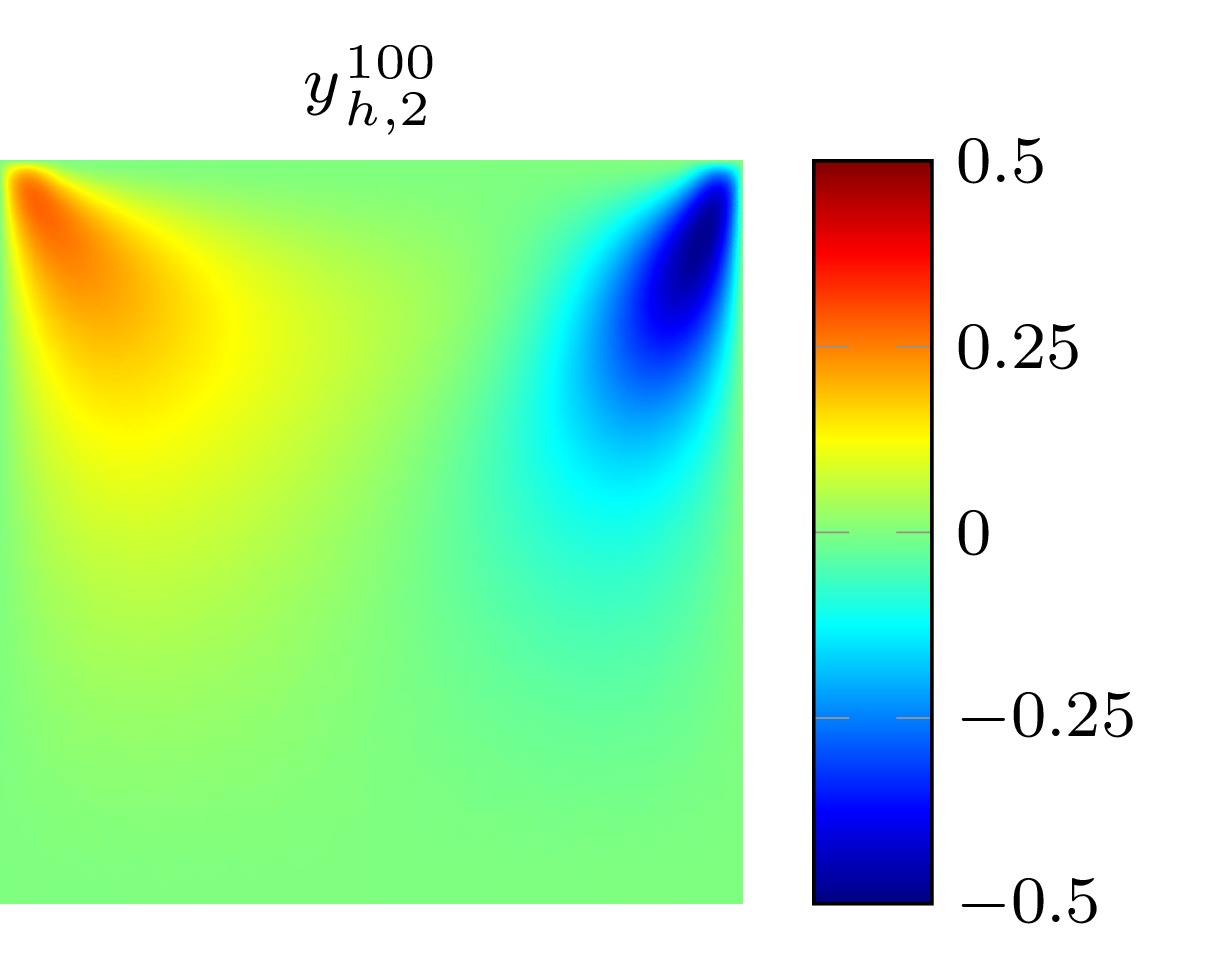}
  \else
    \begin{tikzpicture}[baseline, trim axis left]
      \begin{MyPictureWithColorbar}[title = {$\velhy[100]$},colorbar,colorbar style = {ytick={-0.5,-0.25,0,0.25,0.5}}]
        \input{main15_v_101.tex}
      \end{MyPictureWithColorbar}
    \end{tikzpicture}
  \fi
\\
  \renewcommand{\minVelY}{-1}
  \renewcommand{\maxVelY}{1}
  \ifUseExternalPngs
    \includegraphics{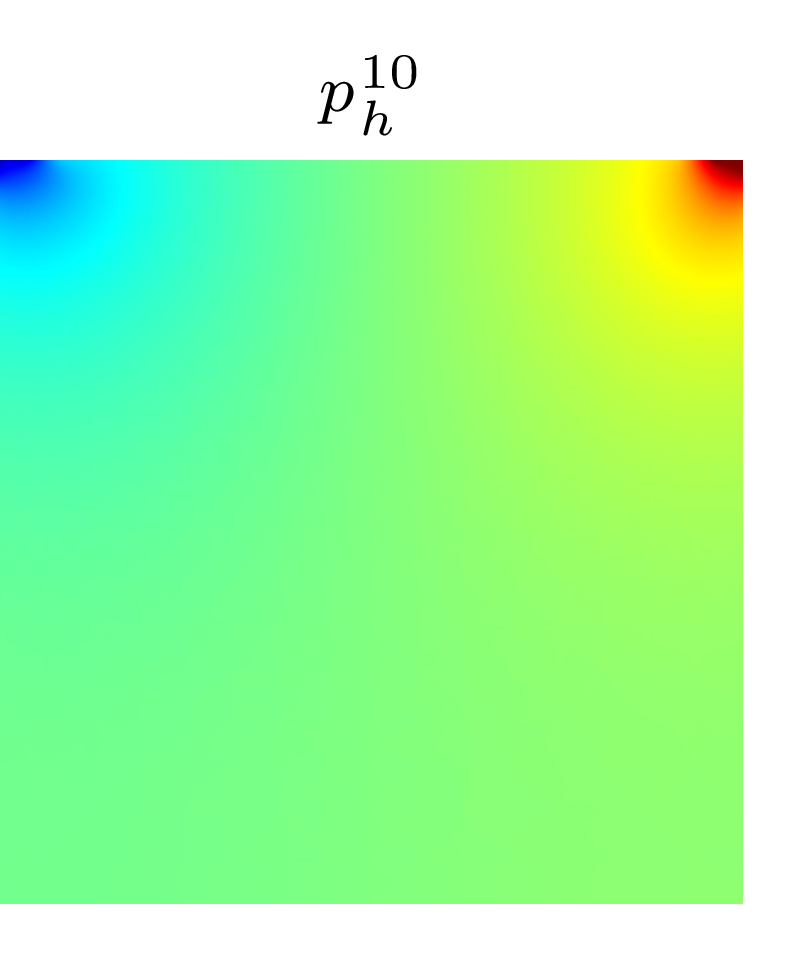}
  \else
    \begin{tikzpicture}[baseline, trim axis left]
      \begin{MyPictureWithoutColorbar}[title = {$\prsh[10]$}]
        \input{main15_p_11.tex}
      \end{MyPictureWithoutColorbar}
    \end{tikzpicture}
  \fi
  \ifUseExternalPngs
    \includegraphics{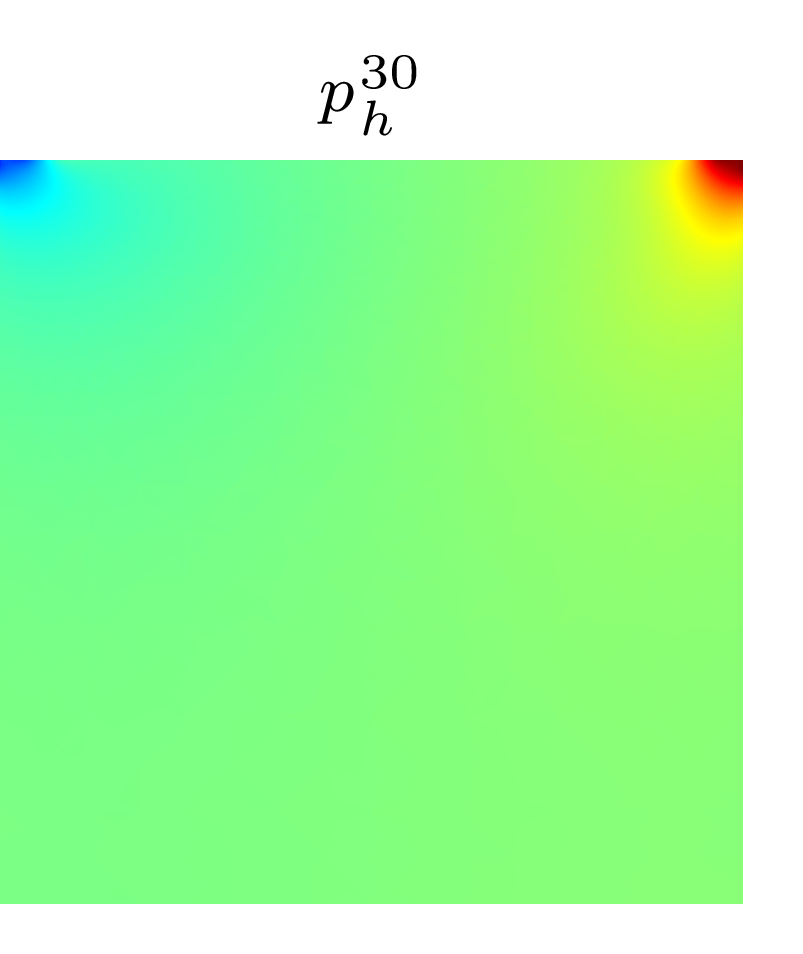}
  \else
    \begin{tikzpicture}[baseline, trim axis left]
      \begin{MyPictureWithoutColorbar}[title = {$\prsh[30]$}]
        \input{main15_p_31.tex}
      \end{MyPictureWithoutColorbar}
    \end{tikzpicture}
  \fi
  \ifUseExternalPngs
    \includegraphics{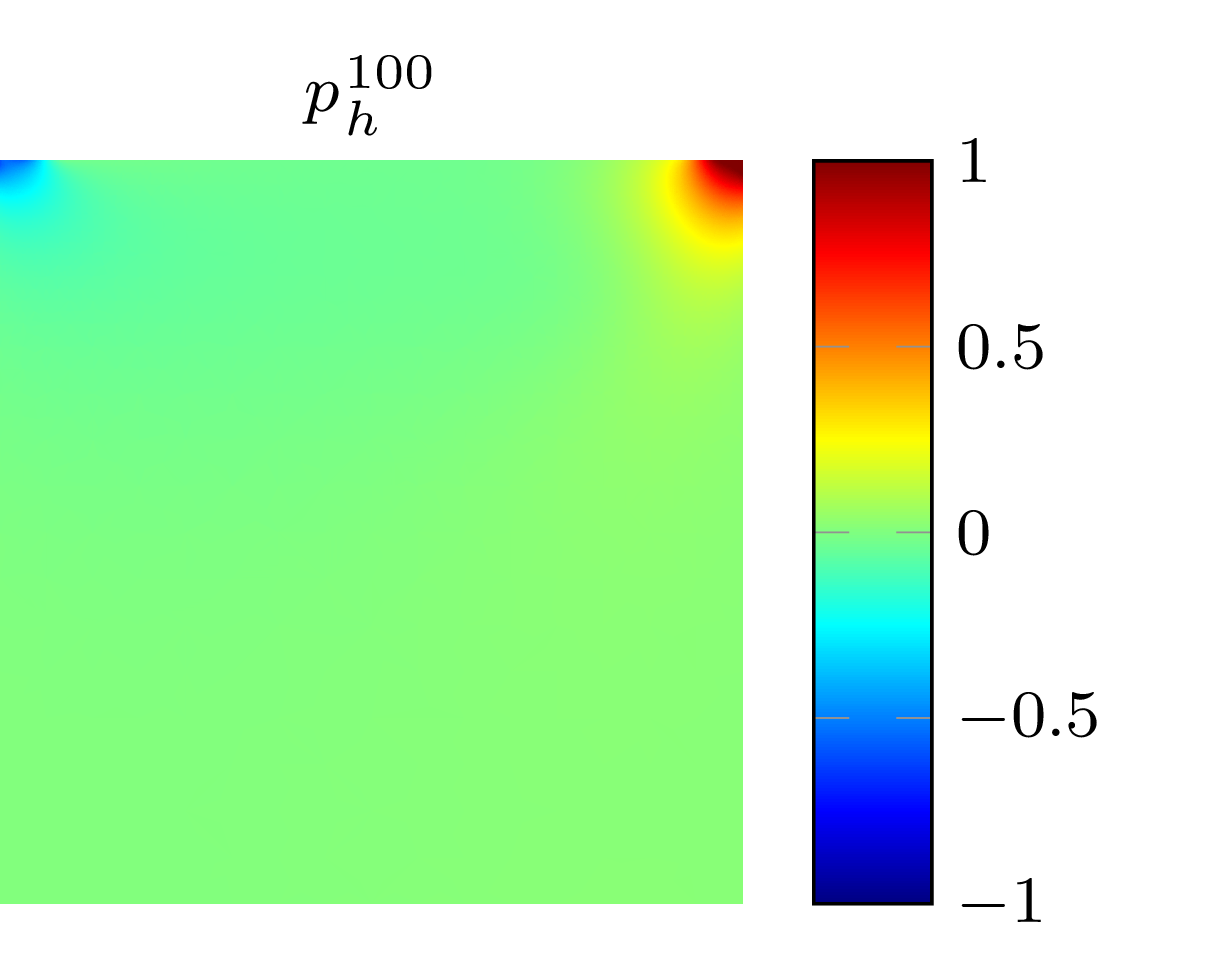}
  \else
    \begin{tikzpicture}[baseline, trim axis left]
      \begin{MyPictureWithColorbar}[title = {$\prsh[100]$},colorbar,colorbar style = {ytick={-1,-0.5,0,0.5,1}}]
        \input{main15_p_101.tex}
      \end{MyPictureWithColorbar}
    \end{tikzpicture}
  \fi
  \caption{Adapted triangulations $\triangulation{j}$, corresponding finite element velocities $\velh[j] = (\velhx[j],\velhy[j])^T$ and pressures $\prsh[j]$ (top to bottom). The columns represent times $t=0.1, 0.3, 1.0$ (left to right).\label{fig:solution}}
\end{center}
\end{figure}

\subsection{Model order reduction}

We need to specify suitable Hilbert spaces and snapshot weights in order to define the construction of POD basis functions from the adaptive finite element snapshots according to \cref{sec:pod}. We also require a set of reference spaces to derive our proposed reduced-order models.

Regarding the choice of Hilbert spaces, we choose the ones that are already used in the weak form and in the finite element error estimator. This means, we take $\hilbSnap=\hilbVel$ for the velocity POD and $\hilbSnap=\hilbPrs$ for the pressure POD in terms of \eqref{minPOD} and we choose $\hilbSnap=\hilbVel$ for the divergence-free projection in \cref{problem:leraySaddlePoint}. 

Regarding the choice of weights, we interpret the sum in the POD minimization problem \eqref{minPOD} as a quadrature of a time integral. A reasonable choice is $\alpha_j=\Delta t$ for $j=1,\dots,\nSnapshots$, which is equivalent to a right-sided rectangle quadrature rule. This complies with the interpretation of the implicit Euler scheme as a discontinuous Galerkin method. 

We define the reference pair of finite element spaces $(\hilbRefVel,\hilbRefPrs)\subset(\hilbVel,\hilbPrs)$ as the pair of finite element space which corresponds to the overlay of all adapted meshes. \Cref{fig:mesh} on the right provides a plot of the overlay of all snapshot meshes of the example simulation. The chosen reference spaces are able to exactly represent all functions in the adapted finite element spaces. Our methods also cover other choices of reference spaces, which enables the decoupling of the POD spatial mesh from the snapshot meshes. However, this can lead to additional interpolation errors depending on the respective resolution.

\subsection{Accuracy}

We compare the considered approaches to model order reduction regarding their accuracy depending on the number of velocity basis functions. In the case of the velocity-pressure model, we set the number of pressure basis functions and the number of stabilizer functions equal to the number of velocity POD basis functions. Descriptions of the tested methods are given in \cref{table:descriptionConvergencePlots}.

\begin{table}
  \caption{Descriptions of reduced-order approximations.\label{table:descriptionConvergencePlots}}
  \begin{tabularx}{\textwidth}{@{}p{0.2\textwidth}X@{}}
    \toprule
    method & description\\
    \midrule
    unstable ROM&velocity-pressure reduced-order model of \cref{sec:supremizers}, but no supremizers\\
    naive ROM&velocity reduced-order model of \cref{sec:approach_i}, but using Lagrange interpolations instead of divergence-free projections of the FE solutions \\
    div-free ROM(1) &velocity reduced-order model of \cref{sec:approach_i}\\
    div-free ROM(2) &velocity reduced-order model of \cref{sec:approach_ii}\\
    div-free POD&optimal approximation of the adaptive FE solutions in terms of the reduced basis of \cref{sec:approach_ii}\\
    stabilized ROM(1)&velocity-pressure reduced-order model of \cref{sec:supremizers}\\
    stabilized ROM(2)&velocity-pressure reduced-order model of \cref{sec:supremizers}, but computed from the pressure snapshots according to \eqref{eq:StabilizerLinearCombination}\\
    stabilized POD&optimal approximation of the adaptive FE solutions in terms of the reduced basis of \cref{sec:supremizers}\\
    \bottomrule  
  \end{tabularx}
\end{table}

\begin{figure}
  \begin{center}
    \ifUseExternalPngs
      \includegraphics{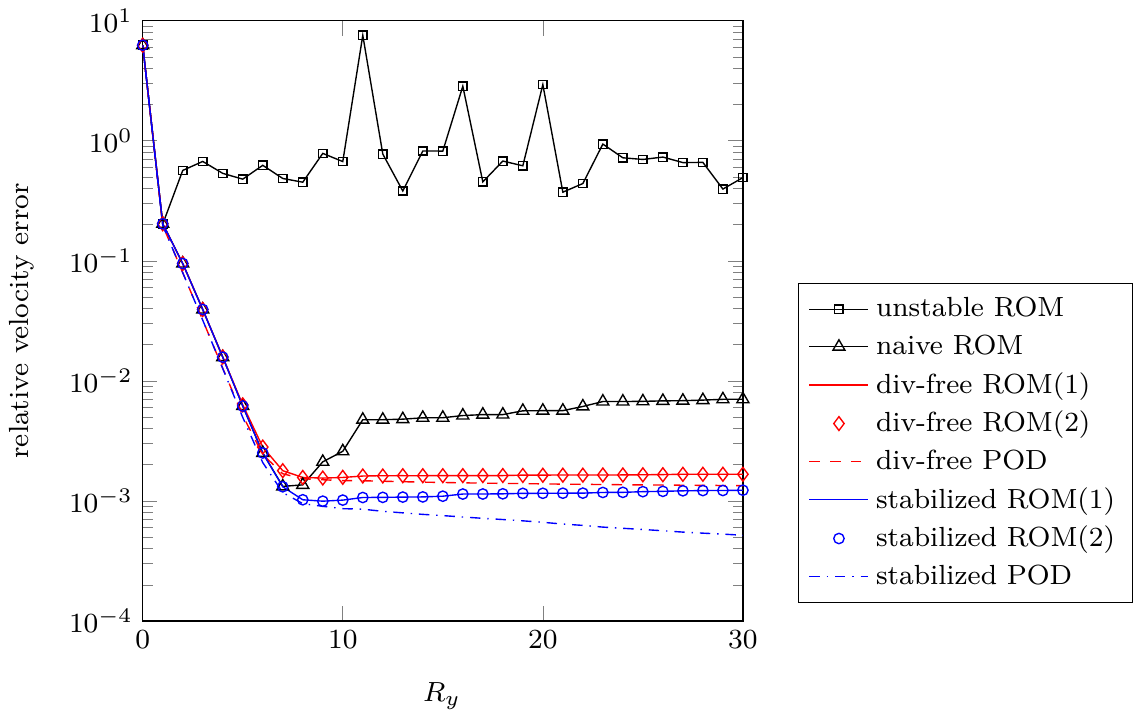}
    \else
      \begin{tikzpicture}
        \begin{semilogyaxis}[
          scale only axis = true,
          width = 0.5\textwidth,
          height = 0.5\textwidth,
          xlabel={$\nPodVel$},
          ylabel={relative velocity error},
          xmin=0, xmax=30,
          ymin=1e-4, ymax=1e1,
          legend entries={
            unstable ROM,
            naive ROM,
            div-free ROM(1),
            div-free ROM(2),
            div-free POD,
            stabilized ROM(1),
            stabilized ROM(2),
            stabilized POD},
          legend style={
            at={(1.65,0.03)},
            anchor=south east,
            cells={anchor=west},
            font=\footnotesize},
          cycle list={
            {black,mark=square, mark options={scale=0.6}},
            {black, mark=triangle},
            {red},
            {red,mark=diamond,only marks},
            {red,dashed},
            {blue},
            {blue, mark=o, mark options={scale=0.75},only marks},
            {blue,dashdotted}}]
          \addplot table [x=R,y=Unstable]       {main15_errorRomAdaptive.dat};
          \addplot table [x=R,y=Naive]          {main15_errorRomAdaptive.dat};
          \addplot table [x=R,y=Leray]          {main15_errorRomAdaptive.dat};
          \addplot table [x=R,y=LerayOfPod]     {main15_errorRomAdaptive.dat};
          \addplot table [x=R,y=LerayOfPod]     {main15_errorProjectionAdaptive.dat};
          \addplot table [x=R,y=Stabilized]     {main15_errorRomAdaptive.dat};
          \addplot table [x=R,y=SupOfSnapshots] {main15_errorRomAdaptive.dat};
          \addplot table [x=R,y=Stabilized]     {main15_errorProjectionAdaptive.dat};
        \end{semilogyaxis}
      \end{tikzpicture}
    \fi
  \end{center}
  \caption{Relative velocity errors of different reduced-order approximations in the sense of \eqref{eq:relErr}, depending on the number of velocity POD basis functions. \label{fig:convergence}}
\end{figure}

We measure the error in the reduced-order velocity approximations with the relative norm implied by the velocity POD, namely
\begin{equation}\label{eq:relErr}
  \text{rel err}_{\vel} = (\sum_{j=1}^\nSnapshots \Delta t \| \velh[j]-\velPOD[j] \|_{\hilbVel}^2)^\frac12/(\sum_{j=1}^\nSnapshots \Delta t \| \velh[j] \|_{\hilbVel}^2)^\frac12.
\end{equation}
The results are provided in \cref{fig:convergence}.

We observe that the relative errors of our proposed approaches show an exponential decay up to $\nPodVel=6$. Thereafter, the decay stagnates at an error slightly above $10^{-3}$. The stagnation of the relative errors in the div-free ROMs and the stabilized ROMs is due to the use of space-adapted snapshots and is related to the finite element discretization error. For more details, we refer to \cite{GH17,URL16}.

The divergence-free velocity reduced-order models and the stabilized velocity-pressure reduced-order models perform similar in terms of accuracy depending on the number of velocity basis functions. Considering the additional degrees of freedom associated with the pressure basis functions and the stabilizing functions, however, the velocity-pressure models are more expensive than the velocity models at the same number of velocity basis.

Both variants of the stabilized velocity-pressure reduced-order model give exactly the same results, in accordance with \eqref{eq:StabilizerLinearCombination}. The difference between the two variants of the divergence-free velocity reduced-order model is in the order of about one percent of the error and, therefore, visually not distinguishable.

Reference curves are given by the projections of the finite element velocity solutions onto the velocity bases used in the reduced-order models. We observe that the errors of the proposed reduced-order models are close to the corresponding optimal errors up to the point where the convergence of the reduced-order model stagnates.

To compare our approaches with less sophisticated methods. A naive reduced-order model is derived by using a non divergence-free velocity basis and neglecting the pressure term and the continuity equation. Nevertheless, the Dirichlet condition is implemented with a weakly divergence-free function, as usual in POD-Galerkin modeling with fixed discretization spaces. The initial convergence of the naive approach is on par with our approach, but if the number of basis functions is increased, the solution starts to diverge. As a second simple alternative, we introduce a velocity-pressure model without stabilizers. Such a model provides a reduced-order solution, but the magnitude of its relative error is of order 1 and, thus, not satisfactory.

\subsection{Cost}

We present computation times of the setup and solution of selected reduced-order models to illustrate the main differences in their computational complexity. We restrict our consideration to a fixed number of velocity basis functions, namely $\nPodVel = 30$. The results are presented in \cref{tab:computationTime}.

Regarding the finite element solution times, we have proved that our adaptive finite element implementation is reasonably efficient by verifying that most of the computation time is spent solving linear systems of equations within the Newton iteration. As initial guess for the Newton iteration, we use the solution at the previous time instance, which is sufficiently accurate to give fast convergence in all our test cases. Considering the total setup times of the proposed methods, we find that they are roughly the same and dominated by the cost of computing the snapshots.

The time to solve the reduced-order model for the time-discrete coefficients of the reduced basis expansion is mainly affected by the setup and solution of the reduced-order linear systems appearing within the nonlinear iteration. The setup costs of the corresponding matrices and right-hand sides are dominated by the third-order convective tensor, whose dimension is equal to the number of velocity unknowns. The solution costs amount to factorizing a dense matrix. Because the velocity models have a smaller number of unknowns, they are significantly more efficient than the velocity-pressure models in terms of solution time. This is also reflected in the measured solution times. The comparison of the full-order simulation with the reduced-order solution gives a speedup factor of 3760 for the div-free method and a factor of 1253 for the stabilized approach. In view of multi-query scenarios like uncertainty quantification or optimal control, where the underlying systems are solved repeatedly, we expect a large gain concerning the computational expenses.

%

\begin{table}
  \caption{Computation times in seconds for selected reduced-order models using $\nPodVel=30$. \label{tab:computationTime}}
  \begin{tabularx}{\textwidth}{@{}p{0.2\textwidth}RRRR@{}}
    \toprule
    & div-free ROM(2)& stabilized ROM(1) &  naive ROM\\
    \midrule
    FE solution& 488.87& 488.87& 488.87\\
    reference FE space & 37.32 & 37.32 &37.32\\
    velocity POD & 0.97 & 0.99 &0.97\\
    pressure POD & -- &0.14  & --\\
    div-free projection & 3.10 & -- &  0.86\\
    supremizers & -- & 0.59 &  --\\
    ROM setup & 1.59& 4.51 &  1.59\\
    ROM solution& 0.13& 0.39 &  0.13\\
\bottomrule
  \end{tabularx}
\end{table}

\section{Conclusions}

In this work, we have extended the framework of POD model order reduction for space-adapted snapshots to incompressible flow problems described by the Navier-Stokes equations. In order to derive a stable POD reduced-order model, two approaches are proposed. 

In the first approach, a velocity reduced-order model is derived by projecting the velocity snapshots or, alternatively, the POD basis functions onto a weakly divergence-free space. In this way, the continuity equation in the POD reduced-order model is fulfilled by construction and the pressure term vanishes. The structural advantage of this approach is that the resulting reduced-order velocity solution is weakly divergence-free with respect to a pressure finite element space regardless of the number of reduced basis functions.

In the second approach, a pair of reduced spaces for the velocity and for the pressure are constructed. Stability is ensured by augmenting the velocity reduced space with pressure supremizer functions. The advantage of this method is that it directly delivers a reduced-order approximation of the pressure field, which is required in many practical applications. The reduced-order velocity, however, is only weakly divergence-free with respect to the pressure reduced space. Moreover, the velocity-pressure reduced-order model is less efficient than the velocity reduced-order model due to the additional degrees of freedom associated with the pressure basis and the supremizers.

Our numerical experiments show that the results of both approaches are very similar in terms of the error between the reduced-order velocity solution and the finite element velocity solution depending on the number of velocity POD basis functions. This implies that the velocity reduced-order model is significantly more efficient than the velocity-pressure reduced-order model in terms of computation time.\\

\textbf{Acknowledgements} Carmen Gr\"a\ss{}le and Michael Hinze gratefully acknowledge the financial support by the DFG through the priority programme SPP 1962. Jens Lang was supported by the German Research Foundation within the collaborative research center TRR154 ``Mathematical Modeling, Simulation and Optimisation Using the Example of Gas Networks'' (DFG-SFB TRR154/2-2018, TP B01). The work of Sebastian Ullmann was supported by the Excellence Initiative of the German federal and state governments and the Graduate School of Computational Engineering at the Technische Universit\"at Darmstadt.

\bibliographystyle{spmpsci}      

\bibliography{literature}   

\end{document}